\documentstyle[amsmath,amssymb,amscd,graphicx,psfrag,epsf,amsthm]{article}

\newcommand{\ou}{{\mathcal {O}}}
\newcommand{\PP}{{\mathbb {P}}}
\newcommand{\A}{{\mathbb {A}}}
\newcommand{\C}{{\mathbb {C}}}
\newcommand{\R}{{\mathbb {R}}}
\newcommand{\X}{{\mathcal {X}}}
\newcommand{\Y}{{\mathcal {Y}}}
\newcommand{\Z}{{\mathcal {Z}}}
\newcommand{\cC}{{\mathcal {C}}}
\newcommand{\cI}{{\mathcal {I}}}
\newcommand{\SRC}{separably rationally connected }

\numberwithin{equation}{section}

\newtheorem{thm}{Theorem}

\newtheorem{lem}[thm]{Lemma}
\newtheorem{cor}[thm]{Corollary}

\newtheorem{prop}[thm]{Proposition}

\newtheorem{conj}[thm]{Conjecture}
   
\theoremstyle{definition}
\newtheorem{defn}[thm]{Definition}

\newtheorem{say}[thm]{}
\newtheorem{exmp}[thm]{Example}

\newtheorem{exrc}[thm]{Exercise}

\newtheorem{rem}[thm]{Remark}  
\newtheorem{rems}[thm]{Remarks}         
\newtheorem{ack}{Acknowledgments}        
   
\newtheorem{warning}[thm]{Warning}  
\newtheorem{defn-thm}[thm]{Definition-Theorem}  
\newtheorem{appl}[thm]{Application}  

\theoremstyle{remark}
\newtheorem{claim}[thm]{Claim}

\renewcommand{\c}[0]{{\mathbb C}}  

\renewcommand{\o}[0]{{\mathcal O}}

\renewcommand{\r}[0]{{\mathbb R}} 
\renewcommand{\a}[0]{{\mathbb A}}

\newcommand{\p}[0]{{\mathbb P}}
\newcommand{\f}[0]{{\mathbb F}}
\newcommand{\q}[0]{{\mathbb Q}}
\newcommand{\map}[0]{\dasharrow}
\newcommand{\qtq}[1]{\quad\mbox{#1}\quad}
\newcommand{\spec}[0]{\operatorname{Spec}}
\newcommand{\pic}[0]{\operatorname{Pic}}
\newcommand{\gal}[0]{\operatorname{Gal}}

\newcommand{\codim}[0]{\operatorname{codim}}

\newcommand{\Hom}[0]{\operatorname{Hom}}

\newcommand{\aut}[0]{\operatorname{Aut}}

\newcommand{\hilb}[0]{\operatorname{Hilb}}

\newcommand{\onto}[0]{\twoheadrightarrow}
\newcommand{\comp}[0]{\operatorname{Comp}}

\def\into{\DOTSB\lhook\joinrel\rightarrow}

\begin{document}

\title{Rational Curves on Varieties}
\author{Carolina Araujo and J\'anos Koll\'ar}
\maketitle

\tableofcontents

\bigskip

The aim of these notes is to give an introduction
to the ideas and techniques of handling rational curves
on varieties. The main emphasis is on varieties
with many rational curves, which are the higher dimensional
analogs of rational curves and surfaces.

Sections  one through five and seven closely follow the lectures
of J.\ Koll\'ar given at the 
summer school ``Higher dimensional varieties and
rational points''. The notes were written up by C.\ Araujo and
substantially edited  later. 
\medskip

\noindent {\bf Terminology.}  We follow the terminology of
\cite{H77}.
Thus a variety is  an integral separated
scheme of finite type over a field (not necessarily algebraically closed).
Somewhat inconsistently, 
by a curve we mean a geometrically reduced (not necessarily irreducible) 
separated $1$-dimensional scheme
of finite type over a field.


\section{Rational Curves on Surfaces}

In this section we study rational curves on surfaces. Our main point 
is that rational surfaces contain plenty of rational curves
(at least when we are over an algebraically closed field). 

In the following example we analyze the possible genera of curves 
lying on smooth surfaces in $\p^3$. This  illustrates our point: while 
surfaces of degree $d\leq 3$ (which are rational) contain many rational 
curves, those of degree $d\geq 5$ tend not to 
contain rational curves at all.

\begin {exmp}[Curves on Surfaces in $\PP^3$] Let $S=S_d$ be 
a smooth surface of degree $d$ in $\PP_{\c}^3$. By the Noether-Lefschetz 
Theorem, if $d\geq4$ and $S_d$ is very general,
any curve $C\subset S_d$ is a complete intersection. This means that there 
is a surface $H_m$ of degree $m$ in $\p^3$, for some $m$, such that 
$C=S_d\cap H_m$ scheme-theoretically.
(See \cite[I]{Lo91} for a survey of different proofs of this 
result.)

Let $\tilde{C}\to C$ be the normalization  of $C$. We ask the following
question: for fixed $d$ and $m$, what are the possible values of 
the genus $g(\tilde{C})$ of $\tilde C$? 

The biggest possible value of $g(\tilde{C})$ is attained when $C=S_d\cap H_m$ 
is already nonsingular. In this case 
$g(\tilde{C})=g(C)=\frac{1}{2}md(m+d-4)+1$.

Lower bounds for $g(\tilde{C})$ are much more interesting and difficult to 
get.

Let $H$ denote the hyperplane class in $\PP^3$, so that $H_m \sim mH$.

By Riemann-Roch we have:
$$
h^0(S,\ou_S(mH)) = m^2 \frac{d}{2} - m \frac{d(d-4)}{2} + \chi(\ou_S).
$$

By the adjunction formula we have:
$$
p_a(C) = m^2 \frac{d}{2} + m \frac{d(d-4)}{2} + 1.
$$

The ``expected'' maximal number of singular points of $C$ is 
$h^0(S,\ou_S(mH)) - 1$.  By ``expected'' we mean that
this is not a theorem but it is frequently correct.
Here is an argument to support it.
 
Requiring that the curve has a singularity 
at a given point imposes 3 conditions. We have a 2-dimensional 
family of points. Thus, requiring that the curve has a singularity at
some point imposes 1 condition.

Putting these formulas together, we obtain that the expected 
smallest genus of $\tilde{C}$ is

$$
md(d-4) + 2 - \chi(\ou_S).
$$

If $d \geq 5$, this grows linearly with $m$. So
we expect that, given $g$, there are only 
finitely many values of $m$ for which $|\ou_S(m)|$ contains a curve
of geometric genus $g$. It is important to keep in mind that this 
``conclusion'' may not hold for every  surface $S_d$. 
However, this is essentially correct for very general surfaces
(although this is not at all easy to prove). In fact, 
by~\cite[Theorem 2.1]{Xu94}, on a very general $S_d\subset \PP_{\c}^3$,
any irreducible 
curve of type $S_d\cap H_m$ has geometric genus $\geq \frac{1}{2}dm(d-5)+2$.

If $d=2$ or $3$, not every curve on $S_d$ is a complete intersection
(for instance, the lines on smooth quadrics and cubics are not complete
intersections). 
One can show that for $d \leq 3$ the linear system 
$|\ou_S(m)|$ contains 
rational curves for every $m$. Instead of proving this fact here, we only
prove that, for $d\leq 3$, $S_d$ contains plenty of rational curves 
(this argument works over any algebraically closed field).

\begin{enumerate}
\item $d=1$. In this case $S_1 \cong \PP^2$, and for 
every $m$ there are plenty of rational curves of degree $m$ (i.e., 
rational curves that can be obtained as $C= S_1 \cap H_m$). These are 
given by maps $\phi: \PP^1 \to \PP^2$ sending $(s:t)$  to  
$(\phi_0(s,t):\phi_1(s,t):\phi_2(s,t))$,
where $\phi_0(s,t)$, $\phi_1(s,t)$ and $\phi_2(s,t)$ are homogeneous 
polynomials of degree $m$.

\item $d=2$. In this case $S_2$ is a smooth quadric surface in 
$\PP^3$, and the stereographic projection from a point in $S_2$ gives 
a birational equivalence $S_2 \sim \PP^2$. 

\item $d=3$. In this case $S_3$ is a smooth cubic 
surface in $\PP^3$, which, again, is birational to $\PP^2$
(it can be realized as $\p^2$ blown up at 6 points). 
\end{enumerate}
 \end{exmp}

\begin{exmp}[Rational curves on $\p^n$]
Given finitely many points $P_1,\dots , P_m$ in $\PP^n$,
let us write down a rational curve 
$C\subset \p^n$ passing through these points.
Set $P_i=(a_{i0}:\dots :a_{in})$. 
Pick any points $c_1,\dots,c_m\in \a^1$.
The Lagrange 
interpolation formula gives  degree $m-1$
polynomials $f_0,\dots,f_n$ such that
$f_j(c_i)=a_{ij}$ for every $i,j$.
These determine a morphism $F:\a^1\to \p^n$
such that $F(c_i)=P_i$.

Note that if the $P_i$ and $c_i$ are defined over a field $k$
then $F$ has coordinate functions defined over $k$.
 
Assume next that we fix not only the points $P_i$
but also tangent directions at each $P_i$.
It is not hard to see that there is a rational curve 
$C\subset \p^n$ through these points 
with the preassigned tangent directions. 

More generally, we can even specify finitely many terms of the
Taylor expansion of  $F:\a^1\to \p^n$
at each point $c_i$.
\end{exmp}

The next example shows that
any smooth rational variety has an even stronger  property.

\begin {exmp}[Constructing curves on rational varieties]
 Let $X$ be a smooth projective rational 
variety of dimension $n$, and  $\Phi:\p^n\map X$ 
 a birational map given by polynomial coordinate functions
$\phi_0,\dots,\phi_N$.

Let $C$ be a smooth projective curve with finitely many points
$c_i\in C$. For each $i$ let us specify
finitely many terms of the Taylor expansion near $c_i$ of a map
from $C$ to $X$. (That is, we choose morphisms
$g_i:\spec \o_{c_i,C}/m_i^{e_i}\to X$, where 
$\o_{c_i,C}$ is the local ring of $C$ at $c_i$ with maximal ideal $m_i$.)

Our aim is to find a morphism $C\to X$ that has the
specified Taylor expansion at the points $c_i$.
Moreover, if the data are defined over a field $k$ we want
$C\to X$ to be defined over $k$.

The first step is to reduce the problem to $X=\p^n$.
In order to do this, we extend each $g_i$ to something bigger.
If we are 
over $\c$, we extend each $g_i$ to a small analytic neighborhood
of $c_i$. Over general fields this corresponds to
extending each $g_i$ to the completion  $\hat{\o}_{c_i,C}$.
The completion of a 1-dimensional smooth local ring is
always the power series ring in one variable, thus these
can be best given by the pull back maps
$$
G_i^*:\o_{P_i,X}\to \hat{\o}_{c_i,C}\cong  k[[t]].
$$
By choosing the extension general enough, we may assume that
the image of $\spec \hat{\o}_{c_i,C}$ is not contained
in the exceptional locus of $\Phi^{-1}$. Thus $G_i$ can be lifted
to a morphism $H_i:\spec k[[t]]\to \p^n$ such that
$G_i=\Phi\circ H_i$. At the function level this means that
$$
G_i^*(u_j)=W_i(t)\cdot H_i^*(\phi_j)\qtq{for every $j$,}
$$
where the $u_j$ are homogeneous coordinates on $X\subset \p^N$,
$G_i^*(u_j)\in k[[t]]$
and  the extra factor $W_i(t)$ comes from the nonuniqueness of
the projective coordinates.  Let $w_i$ be the order of poles of $W_i$.

Assume now that we choose a different map $H'_i:\spec k[[t]]\to \p^n$
that agrees with $H_i$ up to order $d_i$. By composing with $\Phi$
we get $G'_i:\spec k[[t]]\to X$ and the above formula implies that
$G'_i$ agrees with $G_i$ up to order $d_i-w_i$. 

This implies that if we solve our original problem
for maps $C\to \p^n$ for every  $e_i$ then
we also obtain a solution for the problem $C\to X$.

Now let us consider the problem for $C\to \p^n$.
Choose a line bundle $L$ on $C$ of degree at least
$2g(C)-1+\sum e_i$. Then the restriction map
$$
H^0(C,L)\to \sum_i H^0(C,L\otimes \o_{c_i,C}/m_i^{e_i})
\cong \sum_i \o_{c_i,C}/m_i^{e_i}
$$
is surjective. Thus for every $j$ there is a section
$\sigma_j\in H^0(C,L)$ such that the
restriction of $\sigma_j$ to $\o_{c_i,C}/m_i^{e_i}$
agrees with $g_i^*(x_j)$. The map $C\to \p^n$ given by all these
$\sigma_j$ has the required properties.

\end{exmp}

We should keep this in mind in section~\ref{main_section}, when
 we define rationally
connected varieties,
in order to characterize the higher dimensional analogs of rational 
curves and rational surfaces.


\section{Deformation of Morphisms I}

In this section 
we are interested in the following problem: suppose
we have a rational curve $C$ on a smooth variety $X$
defined over some base field $k$. How can we produce
more rational curves out of $C$? 

More generally, let $X$ be a smooth quasi-projective variety, $Y$ 
a projective variety, and $f_0:Y \to X$ a morphism. We would like to
deform the morphism $f_0$ in order to produce more morphisms from $Y$ 
to $X$. 

That is, we would like to find a map
$$
F:Y \times D \to X \qtq{such that} F|_{Y\times \{0\}}=f_0 \ ,
$$
where $(D,0)$ is a smooth pointed algebraic curve. 
For $t\in D$ we set 
$f_t = F|_{Y\times \{t\}}:Y\to X$, and we also require that
$F$ be nontrivial, that is, $f_t \neq f_0$ (as morphisms) for 
general $t$. 

It turns out that the existence of $F$ depends on the
choice of $D$ in a very unpredictable manner.
To avoid this problem,
we need to replace $D$ with a ``very small'' neighborhood
of $0$ in $D$. In general, no affine neighborhood is small enough.
Technically, the best is to work with 
the spectrum of the formal  power series ring $\spec k[[t]]$,
but this is not very intuitive.

Assume for simplicity that we are over $\C$.
Then we get a geometrically  appealing picture
by working with $D=\Delta$, where $\Delta$ denotes the unit disc in $\C$. 
So we want to find a map
$$
F:Y \times \Delta \to X \qtq{such that} F|_{Y\times \{0\}}=f_0.
$$
Here of course $F$ is complex analytic. We also allow
$\Delta$ to be replaced by a smaller disc at any time.
(All discs in $\C$ are biholomorphic to each other,
so this is only a matter of convenience.)



\begin{figure}[hbtp]
\centering

\psfrag{p}{{\fontsize{20}{20}$p$}}
\psfrag{P}{{\fontsize{20}{20}$\PP ^1$}}
\psfrag{F}{{\fontsize{20}{20}$f_t$}}
\psfrag{f}{{\fontsize{20}{20}$f_0$}}
\psfrag{v}{{\fontsize{20}{20}$v_p$}}
\psfrag{q}{{\fontsize{20}{20}$f_0(p)$}}
\psfrag{X}{{\fontsize{20}{20}$X$}}

\scalebox{0.45}{\includegraphics{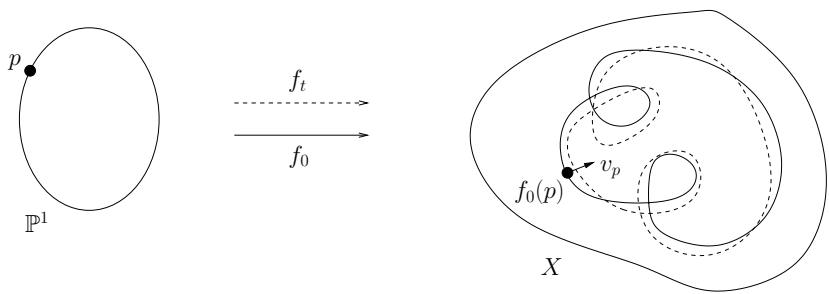}}
\caption{Deforming a rational curve}
\label{deformation}
\end{figure} 


Suppose we have such a map. For each $p\in Y$, $\{f_t(p)\}_{t \in 
\Delta} $ is an analytic arc on $X$. Let $v_p \in T_{f_0(p)}X$ be the
derivative of this arc at the point $f_0(p) \in X$ (see 
Figure~\ref{deformation}). We ``expect'' that this correspondence 
$$
p\in Y \mapsto v_p \in T_{f_0(p)}X
$$ 
 determines a nonzero element 
$\partial F/\partial t\in H^0(Y,f_0^*T_X)$. 
(This is again an expectation only. Given any
$F(y,t):Y\times \Delta\to X$ we can replace it by
 $$
G(y,s):=F(y,s^2):Y\times \Delta\to X
$$ 
to get another morphism  with $\partial G/\partial s=0$.)
In this case, however, the expectation is completely correct,
and we get the following criterion.
 
\begin{center}
{\it{If $H^0(Y,f_0^*T_X)=0$, then $f_0$ cannot be deformed 
in a nontrivial way.}}
\end{center}

\begin {say}[The space of morphisms $Y\to X$] Let $X$ and $Y$ be 
varieties over some 
base field $k$. Suppose we have  
a morphism $f_0:Y\to X$, and we want to study deformations of $f_0$.
In this situation, it is very useful to consider a space 
parametrizing all morphisms from $Y$ to $X$.

For instance, when $Y=X=\p^1$, we would like to say
that isomorphisms $\p^1\to \p^1$ are parametrized by
the group $PGL_2$. More precisely, by this 
we mean that for any field $k$, the points of $PGL_2(k)$
are in a ``natural'' one-to-one correspondence with
the $k$-isomorphisms $\p^1\to \p^1$.

Similarly, when trying to understand morphisms $Y\to X$,
we are looking for a scheme $\Hom(Y,X)$, 
together with a morphism, called {\it universal morphism},
$$
u:\Hom(Y,X)\times Y\to X,
$$
both defined over $k$,
such that for every field extension $K/ k$ and 
every  $K$-morphism $f:Y_K\to X_K$ there is a unique
point $p\in \Hom(Y,X)(K)$ such that
$f=u|_{\{p\}\times Y}$.  The point $p$ is then
denoted by $[f]$.

It turns out that when $X$ is quasi-projective and $Y$ is projective
there is such a  scheme $\Hom(Y,X)$ (see \cite{grot} or \cite[I.1.10]{Ko96}).

It is also very useful to consider morphisms $Y \to X$ sending
a specified point of $Y$ to a specified 
point of $X$. More generally, given $B$ a
closed subscheme of $Y$, and $g:B \to X$ a morphism, we want 
to study morphisms $Y \to X$ extending $g$. Under the same assumptions as 
before, there exists a scheme parametrizing such morphisms. We denote it by 
$\Hom(Y,X,g)$. It is clearly a subscheme of $\Hom(Y,X)$. \end{say}

The scheme $\Hom(Y,X)$ is not always a variety. Usually
it has countably many components. However, each irreducible
component of 
$\Hom(Y,X)$ is of finite type over $k$.

These spaces are very complicated in general. However, if
$H^1(Y,f^*T_X)=0$, $\Hom(Y,X)$ has a very simple description
near the point $[f]$. This is the content of the 
second part of the next theorem.

\begin{thm} \label{ksm} Let $Y$ and $X$ be schemes  over a field $k$,
$Y$ projective and without embedded points 
and $X$ smooth and quasi-projective.
Let $f:Y\to X$ be a morphism. Then
   \enumerate
   \item $T_{[f]}\Hom(Y,X) \cong 
   H^0(Y,f^*T_X)$. 

   \item (Kodaira-Spencer)  If $H^1(Y,f^*T_X)=0$ 
   then $\Hom(Y,X)$ is smooth at $[f]$ and has
      dimension $h^0(Y,f^*T_X)$.

   \item (Mumford) In general, $\dim _{[f]}\Hom(Y,X) \geq 
   h^0(Y,f^*T_X)-h^1(Y,f^*T_X)$. 
\end{thm}

\noindent {\it Proof.} Statements (1) and (3) are special cases of 
\cite[I.2.16]{Ko96}, statements (2.16.1) and (2.16.2) respectively.

Now assume that $H^1(Y,f^*T_X)=0$. It follows from (1) and (3) above that
$\dim _{[f]}\Hom(Y,X)=h^0(Y,f^*T_X)$. Therefore, by \cite[I.2.17]{Ko96},
$\Hom(Y,X)$ is a local complete intersection. Since 
$\dim _{[f]}\Hom(Y,X)=\dim T_{[f]}\Hom(Y,X)$, we conclude that 
$\Hom(Y,X)$ is smooth at $[f]$.
\qed

\medskip

The theorem says that when $H^1(Y,f^*T_X)=0$, there is an irreducible curve
through $[f]$ on $\Hom(Y,X)$ in every direction of
$T_{[f]}\Hom(Y,X)$. This means that any global section of $f^*T_X$ comes from 
a deformation $F:Y\times \Delta \to X$.

This is not the case if $\Hom(Y,X)$ is not smooth at $[f]$.


\section{Deforming Rational Curves}

Now we  apply the theory presented in the previous section to study 
deformations of curves on smooth varieties. We are mainly interested
in deforming rational curves.

Like before, let $X$ be a smooth quasi-projective variety over a field $k$,
but now assume that $Y=C$ is a proper curve over $k$. A crucial
simplification in this case
is that $h^0(C,f^*T_X)-h^1(C,f^*T_X)$ is just
the Euler characteristic $\chi(C,f^*T_X)$, which can be computed by
Riemann-Roch:
$$
\chi(C,f^*T_X) = \deg(f^*T_X) + \chi(\o _C)\dim X = 
-(f_*[C]\cdot K_X) + (1-p_a(C))\dim X.
$$
Hence we obtain that
$$
\dim_{[f]}\Hom(C,X) \geq -(f_*[C]\cdot K_X) 
+ (1-p_a(C))\dim X.
$$

\begin{exmp}  \label{def_fano}
A smooth projective variety $X$ is called \emph{Fano} 
if its anti-canonical class $-K_X$ is ample (for instance, we can take for
$X$ a smooth hypersurface in $\PP^n$ of degree $d\leq n$). 
In this case $(C\cdot K_X)<0$ for every 
curve $C\subset X$, and we expect curves on $X$ to move a lot.
For instance, if $f:\PP^1\to X$ is a nonconstant morphism,
then
$$
\dim_{[f]}\Hom(\PP^1,X) \geq -(f_*[\PP^1]\cdot K_X) 
+ \dim X\geq 1+\dim X.
$$
The automorphisms of $\PP^1$ account for a 3-dimensional family
of maps. Thus, if $\dim X\geq 3$, every rational curve on $X$
moves. (We really need $\dim X \geq 3$ here. Lines on a smooth cubic surface  
in $\p^3$ provide an example of nonmoving rational curves on a Fano variety
of dimension $2$.)
In general, the larger the degree of the curve, the more it moves.

When $X$ is a Fano variety defined over an algebraically closed field, 
a result of Mori shows that 
 through any point of $X$ there passes a rational curve 
(the next exercise shows this in the hypersurface case, when $d\leq n-1$).
\end{exmp}

\begin{exrc} Let $X_d \subset \p^n$ be a smooth hypersurface 
of degree $d\leq n-1$ over an algebraically closed field. 
Show that through every point of $X_d$ there is 
an $(n-1-d)$-dimensional family of  lines.
Compute the dimension of the family of all lines on $X_d$
and compare it with the above estimate.
\end{exrc}

Let us   study the case $C\cong \PP^1$ in more detail. 
Every vector bundle on $\p^1$ splits into a direct sum of line bundles
(see~\cite[V. Exercise 2.6]{H77}). Therefore we can write, 
for suitable integers $a_i$,
$$
f^*T_X \cong \sum_{i=1}^{\dim X}\ou_{\PP^1}(a_i).
$$

\begin{defn} We say that $f:\PP^1\to X$ is \emph{free} if all the $a_i$
are $\geq 0$. 

We say that $f:\p^1\to X$ is \emph{very free} if all the $a_i$ are $\geq 1$.
\end{defn}

\begin{rem}\label{free.moves.rem}
Intuitively, a free curve can
be moved in every direction, while a very free curve still moves freely
in any direction with one point fixed. 

In particular, if $Z\subset X$ has codimension $\geq 2$,
any free curve  $f:\PP^1\to X$
 has nearby deformations that are disjoint from $Z$.
If $f:\PP^1\to X$
is very free, then it
 has nearby deformations $f_t$ such that
$f_t(0)=f(0)$ and $f_t(\PP^1\setminus\{0\})$ 
is disjoint from $Z$.
Moreover, the tangent spaces of the images $f_t(\PP^1)$ at $f(0)$
sweep out an open subset of the tangent space of $X$ at $f(0)$.

These results are easy to obtain from
Propositions \ref{smoothpointisfree} and
\ref{smoothpointisfree.vfree}. See \cite[II.3]{Ko96}
for detailed proofs.
\end{rem}

Notice that a morphism $f:\p^1\to X$ is free (resp.\ very free)
if and only if 
$H^1(\PP^1,f^*T_X(-1))=0$ (resp.\ $H^1(\PP^1,f^*T_X(-2))=0$). By the 
Semicontinuity Theorem (see~\cite[III.12.8]{H77}), nearby deformations
of a free (resp.\ very free) rational curve are still free
(resp.\ very free).

\begin{prop} \label{smoothpointisfree} Let $X$ be an $n$-dimensional smooth 
quasi-projective variety over an algebraically closed field $k$. 
Let $f:\PP^1\to X$ be a nonconstant morphism.
   \begin {enumerate}
   \item If $f$ is free, there exists a variety $Y$ of 
   dimension $n-1$, a point $y_0\in Y$, and a dominant morphism 
   $F:\PP^1 \times Y\to X$ such that $F|_{\PP^1\times \{y_0\}}=f$.
   \item Assume there exists a variety $Y$ and a dominant morphism
   $F:\PP^1 \times Y\to X$. If $k$ has
   characteristic $0$, 
   $f_y:=F|_{\PP^1\times \{y\}}$ is free for general $y\in Y$. 
   \end{enumerate}
\end{prop}

\noindent {\it Proof.} 
   (1) Let $Z$ be a suitable neighborhood of $[f]$
   in $\Hom(\PP^1,X)$. Since $f$ is free,
   any given point of $f(\PP^1)$ can be moved in arbitrary direction.
   This implies that the evaluation morphism $F:\PP^1\times Z\to X$ is 
     submersive at every point of  $\PP^1 \times \{[f]\}$, thus it
          is also dominant. (See \cite[II.3.5.3]{Ko96} for a rigorous proof.)
   If $Y\subset Z$ is a general subvariety of dimension $n-1$ 
   then  $F|_{\PP^1\times Y}:\PP^1\times Y\to X$ is still dominant.
   \medskip 

   \noindent (2) Since $F$ is dominant, there is a dense open subset 
   $U\subset \p^1\times Y$ such that $F:U\to X$ is smooth (here we use
   the characteristic 0 assumption). Therefore, the induced map
   $dF: T_{\PP^1\times Y}\to F^*T_X$ is surjective on $U$.

   Let $\pi_2:\p^1\times Y\to Y$ denote the second projection, 
   and assume that $y \in \pi_2(U)$. Then 
   $$
   \ou_{\PP^1}(2) \oplus \ou_{\PP^1}^{\dim Y}\cong 
   T_{\p^1\times Y}|_{\PP^1\times\{y\}} 
   \to f_y^*T_X \cong \sum_{i=1}^n \ou_{\PP^1}(a_i)
   $$
   is surjective on an open set. Therefore all the $a_i$ are $\geq 0$. \qed

\begin{cor} \label{xfree} Let $X$ be a smooth quasi-projective variety 
over an algebraically closed field of characteristic 0. 
Let $x\in X$ be a very general point, i.e., $x$ is outside a countable union
of proper closed subvarieties of $X$. Then 
any morphism $f:\PP^1\to X$ whose image passes through $x$ is free.\end{cor}

\noindent {\it Proof.} For each (of countably many) irreducible component 
   $Z$ of $\Hom(\PP^1,X)$, consider the evaluation 
   morphism $F_Z:\PP^1\times Z\to X$. 
   Since being free is an open condition, the points of $Z$ parametrizing
   nonfree morphisms form a closed subscheme $Z'\subset Z$ (notice
   that we may have $Z'=Z$). By Proposition~\ref{smoothpointisfree},
   $F_Z|_{\p^1\times Z'}:\p^1\times Z' \to X$ is not dominant, and thus
   $\overline{F_Z(\p^1 \times Z')}\neq X$.

   Set $V_Z=\overline{F_Z(\p^1 \times Z')}\subsetneq X$.
   Then any morphism $f:\PP^1\to X$ whose image passes through a point 
   in $X\setminus \bigcup\limits_ZV_Z$ is free. \qed

\begin{rem}\label{countable.rem}
 This corollary says little  when $X$ is a variety over a countable 
field $k$. Indeed, in this case we may have $X(k)=\bigcup\limits_ZV_Z(k)$
(notation as in the proof above).
No such example is known, and this leads to the  following open problem.

Is there a smooth projective variety $X$ over $\bar{\q}$ such that
through every point of $X$ there  is a nonfree rational curve?

A related question is the following.

Is there a smooth projective variety $X$ over $\bar{\q}$ such that
through every point of $X$ there  is a rational curve,
but rational curves do not cover $X(\c)$?

Variants of the 
 latter are also interesting, and possibly easier, over finite fields.
 \end{rem}

We have an analogue of Proposition~\ref{smoothpointisfree} 
for very free morphisms. The proof is very similar to the proof of 
Proposition~\ref{smoothpointisfree}, and we leave it
to the reader.

\begin{prop}\label{smoothpointisfree.vfree}
Let $X$ be a $n$-dimensional smooth 
quasi-projective variety over an algebraically closed field $k$. 
Fix $x_0\in X$. Let $f:\PP^1\to X$ be a nonconstant morphism such that
$f(0)=x_0$.
   \begin {enumerate}
   \item If $f$ is very free, there exists a variety $Y$ of 
   dimension $n-1$, a point $y_0\in Y$, and a dominant morphism 
   $F:\PP^1 \times Y\to X$ such that $F|_{\PP^1\times \{y_0\}}=f$
   and $F|_{\{0\}\times Y}\equiv x_0$.
   \item Assume there exists a variety $Y$ and a dominant morphism
   $F:\PP^1 \times Y\to X$ such that $F|_{\{0\}\times Y}\equiv x_0$. If
   $k$ has characteristic $0$  then 
   $f_y:=F|_{\PP^1\times \{y\}}$ is very free for general $y\in Y$. 
   \qed
   \end{enumerate} 
\end{prop}

\begin{exmp} \label{fano.exmp}
Let $X$ be a smooth Fano variety over a 
field $k$. We mentioned in Example~\ref{def_fano}
that through any point of $X_{\bar k}$ there passes a rational curve,
defined over $\bar k$. 
When $k$ has characteristic $0$, Theorem~\ref{smoothpointisfree}.2 
implies that through a very general point of $X_{\bar k}$ there passes
a free rational curve. In fact, there is a free rational curve through 
every point, as we see in
Theorem~\ref{RCV}.
On the other hand, 
 so far we know nothing about the existence of  rational 
curves on $X$  that are defined over $k$.  In general this is a 
very difficult
problem that  depends on the arithmetic nature of the field $k$.
In the simplest case, when $k=\r$, we 
  consider this problem in 
 Application~\ref{Rvar}. The more general case of local fields
is studied in Section 9.
\end{exmp}


\section{Deformation of Morphisms II}

Sometimes it is useful to consider a wider class of deformations than 
the one considered in the last two sections.
Let $f:Y\to X$ be a morphism from a projective variety $Y$ to a smooth
quasi-projective variety $X$. Now,
instead of changing only the morphism $f$, we 
also want to change $Y$ and $X$ in a predetermined way.

In order to get this setup,
let $(R,0)$ be an irreducible
pointed scheme, and $\X$
and $\Y$ schemes over $R$, with $\X\to R$ smooth and $\Y\to R$ proper
and flat. Let $\X_0$ and $\Y_0$ be the fibers over $0\in R$. 
Starting with a given 
morphism $f_0:\Y_0 \to \X_0$, we would like to obtain deformations
$f_r: \Y_r \to \X_r$.

Like before, there exists an $R$-scheme $\Hom(\Y /R,\X /R)$ parametrizing
morphisms $\Y \to \X$. This means that there is a
universal morphism
$$
u:\Hom(\Y /R,\X /R)\times_R\Y\to \X
$$
such that for every $R$-scheme $S\to R$ and for every
$S$-morphism
$f:S\times_R\Y\to S\times_R\X$ 
there is a unique $R$-morphism $[f]:S\to \Hom(\Y /R,\X /R)$
such that $f$ is the composition
$$
f:S\times_R\Y\stackrel{(id_S,[f],id_{\Y})}{\longrightarrow}
S\times_R\Hom(\Y /R,\X /R)\times_R\Y
\stackrel{(id_S,u)}{\longrightarrow} S\times_R\X.
$$

The fiber of the map 
$\Hom(\Y /R,\X /R) \to R$ over the point $0\in R$ is precisely
$\Hom(\Y_0,\X_0)$:
$$
   \begin{CD}
   \Hom(\Y_0,\X_0) @>>> \Hom(\Y/R,\X/R)  \\
   @VVV @VVV \\
   \spec k(0) @>>> R
   \end{CD}
$$
(Here $k(0)$ denotes the residue field of $0$ on R.)

The next result is an infinitesimal description of the space
$\Hom(\Y /R,\X /R)$ near a point $[f_0]$ for which 
$H^1(\Y_0,f_0^*T_{\X_0})=0$. (For a proof see \cite[I.2.17]{Ko96}
and the proof of Theorem~\ref{ksm}.)

\begin{thm} \label{ksmR} 
Let $(R,0)$ be an irreducible pointed scheme and $\X$, $\Y$ schemes
over $R$, with $\X \to R$ smooth and $\Y \to R$ proper and flat. 
Assume that $\Y_0$ has no embedded points and let 
$f_0:\Y_0 \to \X_0$ be a morphism. Then
   \begin{enumerate}
   \item If $H^1(\Y_0,f_0^*T_{\X_0})=0$ then 
   $\Hom(\Y/R,\X/R)\to R$ is smooth in a neighborhood of $[f_0]$, and all
   fibers in this neighborhood have dimension $h^0(\Y_0,f_0^*T_{\X_0})$.
   \item In general, 
   $$\dim _{[f_0]}\Hom (\Y/R,\X/R)\geq 
   h^0(\Y_0,f_0^*T_{\X_0})-h^1(\Y_0,f_0^*T_{\X_0})+\dim R.\qed
    $$ 
   \end{enumerate}
   \end{thm}

Here is an example of how this theorem can be applied. 
Let $X$ be a smooth variety over a field $k$.
Suppose 
we have a connected union of
rational curves $D \subset X$, and we want to deform it into a single
rational curve. We represent $D$ as the image of a tree of 
rational curves $C$ under a morphism $f:C \to X$. Since $C$ itself is 
reducible, $\Hom(C,X)$ does not give us any new integral curve on $X$.
So we look for an irreducible $k$-scheme $R$ and an $R$-scheme $\Y$  
(proper and flat over $R$) such that the general fiber $\Y_r$  is
an irreducible  rational curve, while the special fiber $\Y_0$ is 
isomorphic to $C$. (We can usually take $R$ to be $\a^1_k$ and $\Y$ to
be the blow up of $\p ^1_k \times \a^1_k$ at suitable points lying
above $0\in \a^1_k$.
A concrete example is worked out in Application~\ref{move.free.appl}.)
 We also set $\X=X \times_k R$. If  
$H^1(C,f^*T_X)=0$, by Theorem~\ref{ksmR}, $\Z= \Hom(\Y/R,\X/R)\to R$ 
is smooth at $[f]$,
and nearby fibers have dimension $h^0(C,f^*T_X)>0$. By choosing a point
in one such fiber, $[f_r] \in \Z_r$, we obtain a deformation of $f$
$$
f_r:\Y_r \to \X_r=X\times \{r\} \ ,
$$
which is a morphism from a rational curve to $X$.

There is one thing we should be careful about: if $k$ is not 
algebraically closed, it is possible 
that nearby fibers $\Z_r$ do not have $k$-points, and this method does
not yield any deformation of $D$ defined over $k$. We  encounter 
this situation in Application~\ref{Rvar}.

\begin{warning}\label{def.over.etale.ext.warn}
It is important to keep in mind that in general the deformations obtained
in this way do \emph{NOT} give maps from $\Y$ to $\X$, not even
 over any open 
subscheme of $R$. Usually, the only thing we can hope to get is that 
$\Hom(\Y/R,\X/R)\to R$ is smooth at $[f]$.
If this is the situation, we claim that there is 
a subscheme $[f]\in T\subset \Hom(\Y/R,\X/R)$ such that
\begin{enumerate}
\item   $T\to R$ is 
quasifinite  and \'etale, mapping $[f]$ to $0$.
\item $f$ extends to a morphism
$$
F_T:\Y\times_RT\to \X\times_RT.
$$

To see this, let 
$t_1,\dots,t_m$ be a regular sequence of generators
for the maximal ideal  in the local ring
of $[f]$ on $\Hom(\Y_0,\X_0)$. (Such a sequence exists
iff $[f]$ is a smooth point.)
By flatness, these can be lifted to $T_1,\dots,T_m$,
in the local ring
of $[f]$ on $\Hom(\Y/R,\X/R)$.
Then $(T_1=\cdots=T_m=0)$ define $T$ in 
a neighborhood of $[f]$.
\end{enumerate}
 Thus, if for instance
$R=\a_{\c}^1$, we  have deformations over an Euclidean
neighborhood of $0$ in $\a_{\c}^1$, realized as an open subset (in the Euclidean
topology) of some $T\to \a_{\c}^1$. However, there may be no Zariski 
neighborhood of $0$ in $\a_{\c}^1$ over which a nontrivial deformation
 exists.
\end{warning}

\begin{appl}[Moving the union of free rational curves]\label{move.free.appl}
Let $X$ be a smooth variety over an algebraically closed field, 
and suppose we have two rational 
curves $D_1$ and $D_2$ on $X$ such that $D_1\cap D_2 \neq \varnothing$. 
When can we deform $D_1\cup D_2$ into a rational curve?

Let $R=\A^1$, $\X=X\times \A^1$, and $\Y$ be the blow up of 
$\PP^1\times \A^1$ at some point on the fiber over $0\in \A^1$. Then 
$C=\Y_0$ is the union of 2 rational curves, $C_1$ and $C_2$, meeting 
transversely in a point $P$, while all other fibers are smooth 
rational curves. We view $D_1\cup D_2$ as the image of a morphism 
$f_0=f_1\vee f_2:C\to \X_0$. (Notice that we do not require that $D_1$
and $D_2$ meet transversely, nor that they intersect in a single point.
If they meet in more than one point, we choose one of them to 
be the image of $P$.)



\begin{figure}[hbtp]
\centering

\psfrag{A}{{\fontsize{20}{20}$R=\A^1$}}
\psfrag{Y}{{\fontsize{20}{20}$\Y$}}
\psfrag{Y0}{{\fontsize{20}{20}$\Y_0$}}
\psfrag{Yr}{{\fontsize{20}{20}$\Y_r$}}
\psfrag{f0}{{\fontsize{20}{20}$f_0$}}
\psfrag{fr}{{\fontsize{20}{20}$f_r$}}
\psfrag{X}{{\fontsize{20}{20}$X$}}
\psfrag{D1}{{\fontsize{20}{20}$D_1$}}
\psfrag{D2}{{\fontsize{20}{20}$D_2$}}

\scalebox{0.45}{\includegraphics{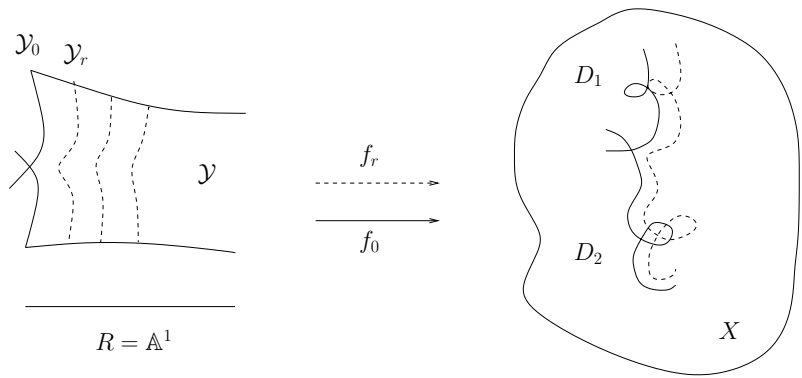}}
\caption{Deforming the union of 2 free rational curves} 
\label{movingunion}
\end{figure} 


By Theorem~\ref{ksmR}, $D_1\cup D_2$ can be deformed into a 
rational curve provided that $H^1(C,(f_1\vee f_2)^*T_{\X_0})=0$. We 
claim that this is the case whenever $f_1$ and $f_2$ are free. 

Indeed, by taking $E=(f_1\vee f_2)^*T_{\X_0}$
in exercise~\ref{h1=0.exrc}, we see that if both $f_1$ 
and $f_2$ are free, then $H^1(C,(f_1\vee f_2)^*T_{\X_0})=0$, and hence 
$D_1\cup D_2$ can be deformed into a rational curve. We actually have more
than that: by taking $E=(f_1\vee f_2)^*T_{\X_0}(-Q)$, where $Q$ is
any point in $C_1\setminus C_2$, we see that  
$H^1(C,(f_1\vee f_2)^*T_{\X_0}(-Q))=0$. By the Semicontinuity Theorem,
$H^1(\p^1,f_r^*T_X(-1))=0$ for nearby deformations $f_r:\p^1 \to X$.
Hence:

\begin{center}{\it{
The union of 2 free rational curves with nonempty intersection  \\ 
can be deformed into a free rational curve.}}
\end{center}
\end{appl}

\begin{exrc} \label{h1=0.exrc}
Let $C=\cup_{i=0}^n C_i$ be a reduced projective curve
of genus zero (i.e., with $h^1(C,\o_C)=0$) and $E$ a vector bundle on $C$.
Assume that
\begin{enumerate}
\item $H^1(C_i,E|_{C_i}(-1))=0$ for $i=1,\dots,n$, and
\item $H^1(C_0,E|_{C_0})=0$.
\end{enumerate}
Then $H^1(C,E)=0$.
\end{exrc}

\begin{appl}[Real Varieties] \label{Rvar}
Let $X$ be a smooth variety defined over $\R$, and assume that $X_{\C}$
contains free rational curves (defined over $\C$). We investigate 
the existence of free rational curves (defined over $\R$) on $X=X_{\R}$.

Assume there exists a point $P\in X(\r)$ and a free rational curve 
$C\subset X_{\c}$ passing through $P$.
If $C$ is defined over $\R$, we are
done. Otherwise, let $\bar{C}$ be its conjugate curve. Then 
$\bar{C}$ 
is also a free rational curve, and their union is defined over $\R$.

Set $R=\A^1_{\R}$, $\X=X_{\R}\times \A^1_{\R}$ and $\Y=(x_0^2+x_1^2=
tx_2^2)\subset \p^2_{\R}\times \a^1_{\R}$. Then $\Y_0$ is the union
of 2 conjugate rational curves meeting transversely in one
point, while all other real 
fibers are smooth rational curves defined over $\R$. The curve 
$\Y_t$ contains $\r$-points for $t>0$ , while  
$\Y_t(\r)=\varnothing$ for $t<0$.
We view 
$C\cup \bar{C}$ as the image of a morphism $f_0:\Y_0 \to \X_0$.

We have $H^1(\Y_0,f_0^*T_{\X_0})=0$. Let $\Z=\Hom(\Y/R,\X/R)$. By 
Theorem~\ref{ksmR},
$\Z\to R$ is smooth at $[f_0]$, and nearby fibers have dimension 
$h^0(\Y_0,f_0^*T_{\X_0})>0$. Now comes a subtle point: it is possible 
that nearby fibers $\Z_r$ have no $\r$-points. We must show that this is not
the case. As in Warning~\ref{def.over.etale.ext.warn},
let  $T\subset \Z$ be a curve through $[f_0]$ 
such that $T\to R$ is smooth  at $[f_0]$.
 We apply the Inverse Function Theorem
for real functions, and conclude 
that $T(\r)\to \r^1$ is a local homeomorphism at $[f_0]$ (in the
Euclidean topology). Thus we obtain many $\r$-points 
in a neighborhood 
of $[f_0]$ in $T$. These points correspond to morphisms 
$f_t:\Y_t \to X\times \{t\}$ that are defined over $\r$. 

Notice that we need that $t>0$ in order to get a curve on $X_{\r}$
containing real points.
The fact that $T\to R$ is \'etale ensures that we can pick such $t$.

By taking $t$ near $0$, we get that $f_t$ is free.
\end{appl}

\begin{rems}
\begin{enumerate}
\item
The argument above requires the existence of a free rational curve
$C\subset X_{\c}$ 
through a real point. 
When $X_{\C}$ contains a free rational 
curve, some dense open subset of $X_{\C}$ is covered by free rational
curves. 
The set of real points, if nonempty, is Zariski dense,
and so we can start the argument.

Problems arise when  $X$ has no real points.
One might suppose that in this case the question 
of existence of free rational curves on $X_{\R}$
does not even make sense.
There can be no real maps $\p^1\to X$ if $X$ has no real points.
There is, however, a quite interesting variant.

Even without real points, $X$ may contain 
empty real conics (i.e., curves obtained from morphisms
$\p_{\r}^2 \supset (x_0^2+x_1^2+x_2^2=0)\to X_{\r}$ defined over $\r$).
If $X(\r)\neq \varnothing$,  
the above deformations corresponding to $t<0$
are such. 
If $X(\r)= \varnothing$, the problem is open.

\item
The argument in Application~\ref{Rvar} 
can be generalized to show that if $X$ is a variety 
defined over a $p$-adic field $k$, and $X_{\bar k}$ contains
a free rational curve passing through a $k$-point $x$, then $X_k$ contains
free rational curves (defined over $k$) through $x$ (see \cite{kol99}).
This is 
because a $p$-adic field also admits an Inverse Function Theorem. 
We return to this in section~\ref{nonACF}.
\end{enumerate}
\end{rems}


\section{Deforming Nonfree Rational Curves I}

In section 3 we studied deformations of free rational curves, and in 
section 4 we saw how to deform unions of free rational curves 
into irreducible ones. Now we are interested in the following problem. 
\begin{center}{\it{
How can we deform a nonfree rational curve on $X$ \\
when $X$ contains many free rational curves?}}
\end{center}

One possible way of dealing with this problem was introduced in 
\cite[1.1.6]{KMM92b}. In this section we explain this method. In the
next section we describe a different approach to the same problem. 
The main trick is to use combs, which we define now.

\begin{figure}[hbtp]
\centering

\psfrag{C0}{{\fontsize{28}{20}$C_0$}}
\psfrag{C1}{{\fontsize{28}{20}$C_1$}}
\psfrag{C2}{{\fontsize{28}{20}$C_2$}}
\psfrag{Cn-1}{{\fontsize{28}{20}$C_{n-1}$}}
\psfrag{Cn}{{\fontsize{28}{20}$C_n$}}
\psfrag{...}{{\fontsize{40}{20}$\dots$}}

\scalebox{0.35}{\includegraphics{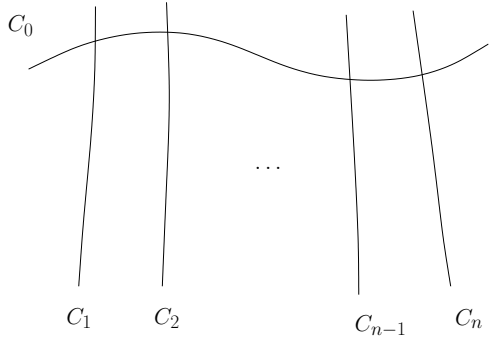}}
\caption{Comb with n teeth}
\label{comb}
\end{figure}

\begin{defn} \label{defcomb}
Let $k$ be an arbitrary field.
A {\it comb} with {\it $n$ teeth} over $k$ is a projective curve
with $n+1$ irreducible components $C_0,C_1,\dots,C_n$ over $\bar k$  
satisfying the following conditions:
\begin{enumerate}
\item The curve $C_0$ is defined over $k$.
\item The union $C_1\cup\dots\cup C_n$ is defined over $k$.
(Each individual curve may not be defined over $k$.)
\item The curves $C_1,\dots,C_n$ are smooth rational curves
disjoint from each other, and each of them meets $C_0$ 
transversely in a single smooth point of $C_0$ 
(which may not be defined over $k$).
\end{enumerate}
The curve $C_0$ is called the {\it handle} of the comb,
and $C_1,\dots,C_n$ are called the {\it teeth}.
A {\it rational comb} is a comb whose handle is a smooth
rational curve.
A comb can be pictured as in Figure~\ref{comb}.
\end{defn}



In section~\ref{nonACF} we  study combs over nonclosed fields
in more detail. For the rest of this section we work over
a fixed algebraically closed field $k=\bar k$ of 
arbitrary characteristic.

We start with an arbitrary  rational curve $D$ on a smooth quasi-projective
variety $X$ over $k$. 
Since $D$ is not free, it is possible that we cannot
move $D$ on $X$. Assume that there are several free rational curves 
$D_1,\dots,D_p\subset X$
meeting $D$ in distinct  
points $d_i\in D$. Then there is a chance that we can deform their
union $D\cup D_1\cup \dots \cup D_p$ into a rational curve. 

Set $R=\A^p$, with affine coordinates $y_i$, and $\X=X\times \A^p$. 
Fix a birational morphism $\p^1\to D$ and 
distinct points $z_1,\dots , z_p\in \p^1$
mapping to $d_1,\dots,d_p$.
Let  $\Y$ be the blow up of $\PP^1\times \A^p$ along
the subvarieties $\{z_i\}\times \{y_i=0\}$, $i=1,\dots ,p$. 
Then $\Y_0$ is a rational comb with handle $C_0$ and
teeth $C_1,\dots ,C_p$. The 
general fiber is a smooth rational curve, but there are also fibers
$\Y_r$, $r\in \a^p$, that are subcombs of 
$C_0\cup C_1\cup \dots \cup C_p$. The number of teeth in this subcomb is
precisely the number of coordinates of $r$ that equal $0$.  
We
view $D\cup D_1\cup \dots \cup D_p$ as the image of a morphism
$f_0:\Y_0 \to \X_0$ mapping $C_0$ birationally
onto $D$ and each $C_i$ birationally onto $D_i$,
$1\leq i \leq p$.

Notice that $\Hom(\Y_0,X)$ only gives us reducible or nonreduced curves 
on $X$. However,
if $\dim_{[f_0]}\Hom(\Y/R,\X/R)>\dim_{[f_0]}\Hom(\Y_0,X)$,
then we can
deform some 
subcurve $D\cup D_{i_1}\cup \dots \cup D_{i_l}$ of $D\cup D_1\cup \dots \cup D_p$
into a rational
curve on $X$. We  show that this is the case if $p$ is large enough.

\begin{rem} \label{limitation}
One limitation of this method is that we only know that 
{\it some} subcurve of $D\cup D_1\cup \dots \cup D_p$ 
 deforms into an irreducible rational curve. In general, we do not even 
know the number of irreducible components of such a subcurve.
\end{rem}

By Theorem~\ref{ksmR} and Riemann-Roch, we have
\begin{equation}\label{5.1.eqn}
\dim_{[f_0]}\Hom(\Y/R,\X/R)\geq -((D \cup D_1 \cup \dots \cup D_p)\cdot K_X)
+\dim X + p.
\end{equation}

Now consider the restriction morphism $\Hom(\Y_0,X)\to \Hom(\PP^1,X)$ 
sending $[h]$ to $[h|_{C_0}]$.
The fiber over $[f]\in \Hom(\p^1,X)$
is isomorphic to 
$$
\prod^p_{i=1} \Hom(\PP^1,X,0\mapsto f(z_i)).
$$ 
By Theorem \ref{ksm}, when
$[f]\in \Hom(\PP^1,X,0\mapsto x)$ is free, 
$$
\dim_{[f]}\Hom(\PP^1,X)=-(f_* [\PP^1]\cdot K_X) + \dim X.
$$
Since we have to impose $\dim X$ conditions in order to fix the image of 
the point $0\in \PP^1$, we get
$$
\dim_{[f]}\Hom(\PP^1,X,0\mapsto x)=\dim_{[f]}\Hom(\PP^1,X)-\dim X=
-(f_* [\PP^1]\cdot K_X).
$$
Therefore we obtain that 
\begin{eqnarray}\label{5.2.eqn}
\dim_{[f_0]}\Hom(\Y_0,X) & \leq & \dim_{[f_0|{C_0}]}\Hom(\PP^1,X)+ \nonumber 
\\ & & -((D\cup D_1 \cup \dots \cup D_p)\cdot K_X) + (D \cdot K_X).
\end{eqnarray}

By inequalities (\ref{5.1.eqn}) and (\ref{5.2.eqn}), we just need to take 
$$
p>\dim_{[f_0|C_0]}\Hom(\PP^1,X) - \dim X + (D \cdot K_X)
$$
in order to have $\dim_{[f_0]}\Hom(\Y/R,\X/R)>\dim_{[f_0]}\Hom(\Y_0,X)$.

We have proved:

\begin{thm} \label{nonfree} Let $X$ be a smooth projective variety
over an algebraically closed field. Let 
$D$ be a rational curve on X (possibly nonfree), and represent it by a 
morphism
$f:\p^1 \to X$. Assume there are $p$ free rational curves 
$D_1,\dots ,D_p$ meeting $D$ in distinct points, with
$p>\dim_{[f]}\Hom(\PP^1,X) - \dim X + (D \cdot K_X)$.
Then there is an integer $l$, $1\leq l \leq p$, and distinct indices
$\{i_1,\dots i_l\}\subset \{1,\dots ,p\}$ such that 
$D\cup D_{i_1}\cup \dots \cup D_{i_l}$ can be deformed into a rational 
curve. \qed \end{thm}

\begin{rem} Given finitely many points on $D$, we can actually
deform some curve $D\cup D_1\cup \dots \cup D_m$ keeping these points
fixed. In this case, however, we may need to attach more free rational
curves to $D$ before we deform it. 

We leave it to the reader to work out the necessary small modifications
when $D$ is not a  rational curve.
\end{rem}

In section~\ref{main_section}, we  see some applications of this theorem.


\section{Deforming Nonfree Rational Curves II}

Let $X$ be a smooth projective variety over an algebraically closed field,
and assume $X$ has very free rational curves. Let $D\subset X$ be a 
 curve. 

In the last section we 
explained the approach of \cite{KMM92b} to the problem of moving 
a curve $D$ on
$X$. First we attach several free rational curves $D_1,\dots ,D_p$
to $D$, and represent $D\cup D_1 \cup \dots \cup D_p$ as the image of
a morphism $f:C\to X$, where $C$ is a comb with $p$ teeth. Then we 
study deformations of $f$.
There are  two limitations of this method.
\begin{enumerate}
  \item As we noted in Remark~\ref{limitation}, by taking $p$ large enough,
  we can only guarantee that some subcurve 
  $D\cup D_{i_1}\cup \dots \cup D_{i_l}$  deforms into an irreducible
  curve. We do not know  which subcurve it is.
  \item If we look at deformations of the morphism $f$, the obstruction
 space   $H^1(C,f^*T_X)$ usually  never vanishes, no matter how large $p$ is.
\end{enumerate}

Recently, \cite{GHS01} introduced a variant to this picture
which is very useful in some applications,
especially in arithmetic questions.

By attaching very free rational curves 
$D_1,\dots ,D_p$ and by looking 
at a modified deformation setting
(in the Hilbert scheme of $X$), we can guarantee that, for $p$ 
large enough, the relevant deformation problem is unobstructed. 
We devote the rest of this section to explaining this new method.

First of all, instead of looking at the space of morphisms to $X$,
we look at the Hilbert scheme of $X$. Below we gather some 
facts about Hilbert schemes. We refer to \cite{grot} and
\cite[I.1 and I.2]{Ko96} for details and proofs.

\begin{say}[Hilbert schemes] \label{hilbert_schemes}
Let $X$ be a projective scheme over a field $k$. There is a $k$-scheme
$\hilb(X)$ parametrizing closed subschemes of $X$.
Like for the Hom schemes, this implies that
 there is a universal subscheme
$U\subset \hilb(X)\times X$, flat over $\hilb(X)$, 
such that 
for every field extension $K/ k$ 
and every
$K$-subscheme $Z$ of $X_K$ there is a unique point
$[Z]\in \hilb(X)(K)$ such that the fiber of $U$ over $[Z]$
is exactly $Z$.

Now assume that $X$ is smooth, and let $Z$ be a closed subscheme of $X$
given by the ideal sheaf $\cI_Z$. 
Assume that $Z$ is a local complete intersection, that is,
$\cI_Z$ is locally generated by $\codim(Z,X)$ elements.
Then $\cI_Z/\cI_Z^2$ is a locally free sheaf on $Z$
and its dual $\Hom (\cI_Z/\cI_Z^2,\o_Z)$ is the 
 {\it normal bundle} $N_Z$ of $Z$ in $X$. We have:
\begin{enumerate}
\item $T_{[Z]}\hilb(X)\cong H^0(Z,N_Z)$.
\item If $H^1(Z,N_Z)=0$, then $\hilb(X)$ is smooth at $[Z]$.
\item In general, $\dim_{[Z]}\hilb(X)\geq h^0(Z,N_Z)-h^1(Z,N_Z)$.
\end{enumerate}
(See \cite{grot} or \cite[I.2.8]{Ko96}.)
\end{say}

Now we go back to our problem. 

Let $X$ be a smooth projective variety
over 
an algebraically closed field
 and $D\subset X$ an irreducible curve.
Our goal is 
to attach  rational curves $D_1,\dots,D_p$ to $D$
so that the resulting comb $D\cup D_1\cup\dots\cup D_p$ deforms 
into an irreducible curve. 

The first problem we face is that $D$ may be quite singular,
in which case we prefer to work with its normalization $\tilde D\to D$.
On the other hand, here it is essential that $D$ be embedded in $X$.
To achieve both, we take an embedding $\tilde D\into \p^3$
and consider the diagonal map $\delta:\tilde D\to X\times \p^3$.
This is an embedding with image $D'\cong \tilde D$
whose projection to the first factor is $D$.
Any moving of $D'$ can be projected to $X$ to give a
deformation of $D$. For most applications this reduction
to the smooth curve $D'$ is satisfactory. There may be, however,
delicate geometric questions for which one needs to be more careful.

So assume $D\subset X$ is a smooth irreducible curve. 
Suppose there are disjoint smooth rational curves $D_1,\dots ,D_p$, 
each meeting $D$ transversely in a single point.
Then we can deform their union $C=D\cup D_1\cup \dots \cup D_p$
into an irreducible curve $C'$ provided that the following conditions hold:
\begin{enumerate}
\item The sheaf $N_C$ is generated by global sections.
\item $H^1(C,N_C)=0$.
\end{enumerate}
Indeed, we have an exact sequence
$$
0\ \to \ N_D \ \to \ N_C|_D \ 
\to \ Q \ \to \ 0 ,
$$
where $Q$ is a torsion sheaf supported at the points $P_i=D\cap D_i$. 
If condition~1 holds, we can find a global section 
$s\in H^0(C,N_C)$ such that, for each $i$, 
the restriction of $s$ to a neighborhood of $P_i$ 
is not in the image of $N_D$. This means that $s$ 
corresponds to a first-order deformation of $C$ that smoothes the nodes
$P_i$ of $C$.
If condition~2 holds, there is no obstruction, and we can find a global
deformation of $C$ that smoothes its nodes $P_i$.

\begin{rem}
In many situations, we would like the resulting irreducible curve 
$C'$ to have some extra properties. For instance, if we start with
a smooth rational curve $D$, we may want $C'$ to be a free 
(resp. very free) curve. In this case, by the Semicontinuity Theorem,
we can just fix a point $P\in D$ (resp. points
$P,Q\in D$) and replace condition (2) above by
\begin{enumerate}
\item[2'.] $H^1(C,N_C(-P))=0$ \ (resp. $H^1(C,N_C(-P-Q))=0$).
\end{enumerate}
\end{rem}

\begin{thm} \label{nonfree2}
Let $X$ be a smooth projective variety 
of dimension at least 3
over an algebraically closed field.
Let   $D\subset X$ be  a smooth irreducible curve and 
 $M$  a  line bundle on $D$. 
Let $C\subset X$ be a very free rational curve intersecting $D$
and let
$\cC$  be a family of rational curves  on $X$ parametrized by
a neighborhood of $[C]$ in $\Hom(\p^1,X)$.

Then there are 
curves $C_1,\dots,C_p\in \cC$ such that  
$D^*=D\cup C_1\cup \dots \cup C_p$ is a comb and
satisfies the following conditions:
\begin{enumerate}
\item The sheaf $N_{D^*}$ is generated by global sections.
\item $H^1(D^*,N_{D^*}\otimes M^*)=0$, where $M^*$ is the unique
line bundle on $D^*$ that extends $M$ and has degree $0$ on the $C_i$.
\end{enumerate}
\end{thm}

We have a lot of freedom in choosing curves from  $\cC$.
Namely, as in Remark~\ref{free.moves.rem},
for a general point $P\in D$
and for a general tangent direction
$\ell\subset T_{P,X}$ there is a curve
$C=C(P,\ell)\in \cC$ passing through $P$ and having tangent direction
$\ell$ at $P$. Moreover, given a 1-dimensional subscheme
$B\subset X$ we can also assume that $C$
intersects $D\cup B$ only at $P$. (This is the only point where we need
$\dim X\geq 3$.)
By looking at the proof, we see that
it is sufficient that there is a curve in $\cC$ having a general 
normal direction $\bar{\ell}\subset N_{P,D}$ at $P$ 
(instead of a general tangent direction $\ell \subset T_{P,X}$). 
This is very useful in the
proof of Theorem~\ref{ghs}.

For arithmetic applications, one has to be even more
careful with the choice of the teeth $C_i$
(see \cite{koll-spec}).
\medskip

\noindent {\it Proof.} 
{\it Step 1.} Here we produce a comb $D'\subset X$ with handle $D$ 
such that $H^1(D,N_{D'}|_D\otimes M)=0$. 

Set $D_0=D$. Suppose we have a comb 
$D_i=D\cup C_1\cup \dots\cup C_i$  
such that $d_i=h^1(D,N_{D_i}|_D\otimes M)\neq 0$. 
We show how to choose a curve $C_{i+1}\in \cC$, disjoint from 
$C_1,\dots,C_i$, forming a comb $D_{i+1}=D_i\cup C_{i+1}$ so that 
$d_{i+1}<d_i$. Then we just repeat this procedure finitely 
many times until we reach $d_p=0$. 
   
By Serre duality, $H^1(D,N_{D_i}|_D\otimes M)$ is dual to
$\Hom(\omega_D^{-1}\otimes M, (\cI_{D_i}/\cI_{D_i}^2)|_D)$,
so we would like to understand how  $(\cI_{D_i}/\cI_{D_i}^2)|_D$
changes with $i$.
 
Pick  any  $C_{i+1}\in \cC$ disjoint from $C_1,\dots, C_i$ and 
intersecting $D$ transversely in a single point $P$.
The tangent direction of $C_{i+1}$ at $P$ corresponds to a normal
direction $l_{C_{i+1}}\subset N_{P,D}$. By explicit computation
we get an 
exact sequence 
$$
0\ \to \ (\cI_{D_{i+1}}/\cI_{D_{i+1}}^2)|_D \ \to \ 
(\cI_{D_i}/\cI_{D_i}^2)|_D \ \to \ l_{C_{i+1}}^{\vee} \ \to \ 0,
$$
where the last map is the composition of the restriction map
$(\cI_{D_i}/\cI_{D_i}^2)|_D \to N_{P,D}^{\vee}$ with the dual map
$N_{P,D}^{\vee}\to l_{C_{i+1}}^{\vee}$.

In particular, we get an  inclusion of groups 
$$
\Hom(\omega_D^{-1}\otimes M, (\cI_{D_{i+1}}/\cI_{D_{i+1}}^2)|_D)\ \into \ 
\Hom(\omega_D^{-1}\otimes M, (\cI_{D_i}/\cI_{D_i}^2)|_D).
$$
So we always have $d_{i+1}\leq d_i$.

Since $d_i\neq 0$, there is a nonzero element 
$\phi \in \Hom(\omega_D^{-1}\otimes M, (\cI_{D_i}/\cI_{D_i}^2)|_D)$.
We are done if we can choose $C_{i+1}$ so that $\phi$ does not come
from an element $\tilde \phi$ of 
$\Hom(\omega_D^{-1}\otimes M, (\cI_{D_{i+1}}/\cI_{D_{i+1}}^2)|_D)$.

For each $P\in D$, let $\phi_P$ denote the composition of $\phi$ with
the restriction map $(\cI_{D_i}/\cI_{D_i}^2)|_D \to N_{P,D}^{\vee}$. 
Since $\phi$ is nonzero, there is a dense open subset $V\subset D$ such that
$\phi_P$ has rank $1$ for every $P\in V$. Let 
$\xi_P\in N_{P,D}^{\vee}$ span the image of $\phi_P$.

Choose a general point 
$P\in V$ not in $C_1\cup \dots \cup C_i$ and a curve $C_{i+1}\in \cC$
through $P$ such that
\begin{enumerate}
\item $C_{i+1}$ intersects $D\cup C_1\cup \dots \cup C_i$
only at $P$, and
\item  The projection of the tangent direction  of $C_{i+1}$ at $P$
to $N_{P,D}$, denoted by $l_{C_{i+1}}$,
is not contained in the kernel of $\xi_P:N_{P,D}\to k$.
\end{enumerate}

It follows from our exact sequence that  
$\phi \in \Hom(\omega_D^{-1}\otimes M, (\cI_{D_i}/\cI_{D_i}^2)|_D)$ 
comes from an element $\tilde \phi$ of
$\Hom(\omega_D^{-1}\otimes M, (\cI_{D_{i+1}}/\cI_{D_{i+1}}^2)|_D)$
if and only if $\xi_P|_{l_{C_{i+1}}}=0$. This does not hold by the 
choice of $C_{i+1}$, and we are done.

Notice that, once we have a comb $D'$ such that
$H^1(D,N_{D'}|_D\otimes M)=0$, any other comb $D^*$
obtained from $D'$ by adding extra teeth  also satisfies 
$H^1(D,N_{D^*}|_D\otimes M)=0$.  

\medskip

\noindent {\it Step 2.} Now we extend the comb $D'$ obtained 
in Step 1 to a comb $D^*$ so that $N_{D^*}|_D$ is 
generated by global sections.

Let $D'$ be the comb obtained in Step 1. Let 
$M'$ be any line bundle on $D$ with degree $-(g+1)$, where $g$ is the 
genus of the curve $D$. We use
the method described above to attach more teeth to $D'$ so that 
the resulting comb $D^*$ also satisfies $H^1(D,N_{D^*}|_D\otimes M')=0$.
By Riemann-Roch, $M'$ embeds in $\ou_D(-P)$ for every $P\in D$. 
So we get that $H^1(D,N_{D^*}|_D(-P))=0$ for every $P\in D$, and this 
implies global generation of $N_{D^*}|_D$.

\medskip

\noindent {\it Step 3.} Finally, we use the fact that the $C_i$ are free
to show that 
$N_{D^*}$ is generated by global sections and that 
$H^1(D^*,N_{D^*}\otimes M^*)=0$

Let $D^*=D\cup C_1\cup \dots \cup C_p$ be the comb obtained in Step 2,
so that $N_{D^*}|_D$ is generated by global sections and 
$H^1(D,N_{D^*}|_D\otimes M)=0$. 

For each $i$ we have an exact sequence
$$
0\ \to \ N_{C_i} \ \to \ N_{D^*}|_{C_i} \ 
\to \ Q_i \ \to \ 0 ,
$$
where $Q_i$ is a torsion sheaf supported at  $P_i=D\cap C_i$.
Since $C_i$ is free, $H^1(C_i,N_{C_i}(-P_i))=0$, and hence 
$H^1(C_i,N_{D^*}|_{C_i}(-P_i))=0$.

We also have the exact sequence 
$$
0\ \to \ \sum_{i=1}^p N_{D^*}|_{C_i}(-P_i) \ \to \ N_{D^*} \ 
\to \ N_{D^*}|_D \ \to \ 0.
$$
Therefore, the vanishing of $H^1(D^*,N_{D^*}\otimes M^*)$ follows from the 
vanishing of $H^1(D,N_{D^*}|_D\otimes M)$ and $H^1(C_i,N_{D^*}|_{C_i}(-P_i))$,
$1\leq i\leq p$.

Similarly, global generation of $N_D^*$ follows from global generation of
$N_{D^*}|_D$ and $N_{C_i}$, $1\leq i\leq p$, plus the vanishing of 
$H^1(C_i,N_{D^*}|_{C_i}(-P_i))$. \qed


\section{Rationally Connected Varieties} \label{main_section}

In this section we introduce the concept of rationally connected varieties.
After some motivation, we present possible definitions of the concept
and show that they are equivalent. Then we study some properties of the class
of rationally connected varieties. This is all done over $\c$ (or any 
uncountable algebraically closed field of characteristic $0$). At the end
of the section we consider subtleties that arise in positive characteristic.

\begin{say}[Motivation] It is rather clear that  rational curves
form a distinguished class among all  algebraic curves. They are the simplest
algebraic curves in many aspects (number theoretical, topological, etc.).

In  dimension two, smooth rational surfaces are the best analogs of 
 rational curves. Their class enjoys several  nice properties:

\begin{enumerate}
   \item (Castelnuovo's Criterion) A smooth projective surface $S$ over
   $\c$ is rational if and only if 
   $H^0(S,\omega_S^{\otimes 2})=H^1(S,\o_S)=0$.
   \item (Deformation Invariance) Let $S \to T$ be a family of
   smooth projective 
   surfaces over $\c$, $T$ connected.
   If $S_0$ is rational, any $S_t$ is rational.
   \item (Noether's Theorem) Let $S$ be a smooth projective surface. If there 
   is a morphism $S \to \p^1$ whose general fibers are rational curves,
   then $S$ itself is rational.
\end{enumerate}

When we move to higher dimensions, and try to classify the ``simplest''
algebraic varieties, things get much more complicated. As a first
try, one could take rational varieties as higher dimensional analogs of 
rational curves and rational surfaces. 
However, there are higher dimensional varieties that behave 
like rational varieties but are not rational. ``Behaving like rational 
varieties'' should be understood as having the relevant properties of
rational varieties.
Smooth cubic 3-folds, for instance, behave like rational varieties in 
many aspects, but they are 
not rational (see \cite[Theorem 13.12]{CG72}). Moreover, the analogs of 
properties 1--3 above all seem to fail for rational varieties 
of dimension 3 and up. Smooth cubic 3-folds satisfy 
$H^i(S,\o_S)=H^0(S,\omega_S^{\otimes m})=0$ for every $i,m\geq 1$. There 
are examples of conic bundles over $\p^2$ which are not rational. It is
conjectured that deformation invariance also fails.

Another possibility would be to consider \emph{unirational varieties}. A
variety of dimension $n$ is unirational if there is a dominant
map $\p^n \map X$. However, the class of unirational varieties
still seems to be too small. It is conjectured that deformation 
invariance fails for the class of unirational varieties.

A new concept was proposed in \cite{KMM92b}, namely the concept of 
\emph{rationally connected varieties}. The definition was based on the fact
that $\p^n$ has many rational curves, and on the expectation that a 
variety $X$ behaves like $\p^n$ if and only if it has plenty of rational curves.
The precise notion is given in Definition-Theorem~\ref{RCV} below.

It seems rather clear by now that in characteristic $0$
rationally connected varieties are
the right higher dimensional analogs of rational curves and rational 
surfaces. This class of varieties enjoys many nice properties: 
\begin{enumerate} \setcounter{enumi}{-1}
\item (Birational invariance) 
    The class of rationally connected varieties is closed under 
   birational equivalence. 
\item (Castelnuovo's Criterion) A smooth projective rationally 
   connected variety $X$ satisfies
\begin{equation*}
H^0(X,(\Omega_X^1)^{\otimes m})=0 \ \ \ for\ \ every\ \  m\geq 1,
\end{equation*}
 and the converse is conjectured to hold.  (This is 
   Theorem~\ref{Castelnuovo}.)
\item (Deformation Invariance)
   Being rationally connected is deformation invariant for smooth 
   projective varieties over. 
   (This is Theorem~\ref{deformation_invariance}.) 
\item (Noether's Theorem)
If $X$ is a smooth variety, and there is a dominant morphism 
   $X \to Y$ to a rationally connected variety $Y$ such that
   the general fibers are rationally connected,
   then $X$ itself is rationally connected. (This is
   Corollary~\ref{rcfibration}.) 
\end{enumerate}

The class of rationally connected varieties contains the class of 
unirational varieties. It should be said, however, that it is still an
open problem to show that these classes are in fact distinct.
That is, there is not a single rationally connected variety over $\c$
that is known not to be unirational. Potential candidates abound,
for instance general quartics in $\PP^4$, any smooth hypersurface of degree $n$
in $\PP^n$ for $n\geq 5$, smooth hypersurfaces of bidegree $(n,2)$
in $\p^k\times \p^2$ for $n\gg k\geq 2$.
\end{say}

\begin{defn-thm} \label{RCV} Let $X$ be a smooth projective variety over
$\C$ (or any uncountable algebraically closed field of characteristic $0$). 
We say that $X$ is \emph{rationally connected} if it satisfies the
following equivalent conditions.

\begin{enumerate}

\item[(1)] There is a dense open set $X^0\subset X$ such that, for 
   every $x_1,x_2\in X^0$, there is a chain of rational curves connecting
   $x_1$ and $x_2$.

\item[(2)] There is a dense open set $X^0\subset X$ such that, for 
   every $x_1,x_2\in X^0$, there is a rational curve through $x_1$ and $x_2$.

\item[(3)] For every $x_1,x_2\in X$, there is a rational curve 
   through $x_1$ and $x_2$.

\item[(4)] For every integer $m>0$, and every $x_1,\dots ,x_m\in X$, 
   there is a rational curve through $x_1,\dots ,x_m$.

\item[(5)] For every integer $m>0$, and every $x_1,\dots ,x_m\in X$, 
   there is a free rational curve through $x_1,\dots ,x_m$.

\item[(6)] There is a very free morphism $f:\PP^1 \to X$, that is, 
   $$
   f^*T_X=\sum_{i=1}^{\dim X} \ou_{\p^1}(a_i) \qtq{with}  a_i\geq 1.
   $$
\end{enumerate}

If $X$ is defined over an arbitrary field $k$ of characteristic $0$,
we say that $X$ is \emph{rationally connected} if $X_K$ satisfies the 
conditions above for some (or equivalently for all) 
uncountable algebraically closed field $K$ extending $k$.

It is possible that the six conditions above are equivalent over
any algebraically closed field of characteristic $0$, 
but this is not known for countable fields, cf.\
Remark \ref{countable.rem}.

\end{defn-thm}

\noindent We of course need to prove that the six versions above are
equivalent. The proof is somewhat more transparent if we introduce
two other variants:

\begin{enumerate}

\item[(1$^+$)] For 
   every $x_1,x_2\in X$, there is a chain of rational curves connecting
   $x_1$ and $x_2$.

\item[(6$^+$)] For   every $x\in X$
 there is a very free morphism $f:\PP^1 \to X$ through $x$.

\end{enumerate}

The implications (5) $\Rightarrow$ (4) $\Rightarrow$ (3) $\Rightarrow$ 
(2) $\Rightarrow$ (1) are obvious,
and so are  (1$^+$) $\Rightarrow$ (1) and (6$^+$) $\Rightarrow$ (6).
\medskip

Assume (1) and let us prove (6).
The argument also shows (1$^+$) $\Rightarrow$ (6$^+$).

Let $X^{\rm free}$ be the set obtained in Corollary~\ref{xfree}, so that any 
rational curve meeting $X^{\rm free}$ is free. 
It is the intersection of countably many dense open subsets.
We repeatedly use two of its properties over $\c$: $X^{\rm free}$ is dense
and for every irreducible subvariety $Z\subset X$, the intersection
$X^{\rm free}\cap Z$ is either empty or dense in $Z$. Thus, if a curve
$C\subset X$ intersects $X^{\rm free}$, a dense set of its deformations
also intersect $X^{\rm free}$.

Pick  a point $x\in X^{\rm free}\cap X^0$ and
let $x'$ be any point in $X^0$.
By hypothesis, there is a chain of rational curves 
$C_1\cup \dots \cup C_n$ connecting $x$ and $x'$. We
show that there is in fact a chain of \emph{free} rational
curves connecting them.

Set $x_0=x$, $x_i=C_i\cap C_{i+1}$ for $0<i<n$, and $x_n=x'$.
We  replace each curve $C_i$ by a free rational curve $C_i'$ passing
through $x_i$ and $x_{i+1}$.

The curve $C_0$ is itself free because it intersects $X^{\rm free}$, so we set 
$C_0'=C_0$.
 
\begin{figure}[hbtp]
\bigskip
\centering

\psfrag{x0}{{\fontsize{28}{20}$x_0$}}
\psfrag{x1}{{\fontsize{28}{20}$x_1$}}
\psfrag{x2}{{\fontsize{28}{20}$x_2$}}
\psfrag{C0}{{\fontsize{28}{20}$C_0$}}
\psfrag{C1}{{\fontsize{28}{20}$C_1$}}
\psfrag{C1'}{{\fontsize{28}{20}$C_1'$}}
\psfrag{D1}{{\fontsize{28}{20}$D_1$}}
\psfrag{Dk}{{\fontsize{28}{20}$D_k$}}
\psfrag{X}{{\fontsize{28}{20}$X$}}

\scalebox{0.35}{\includegraphics{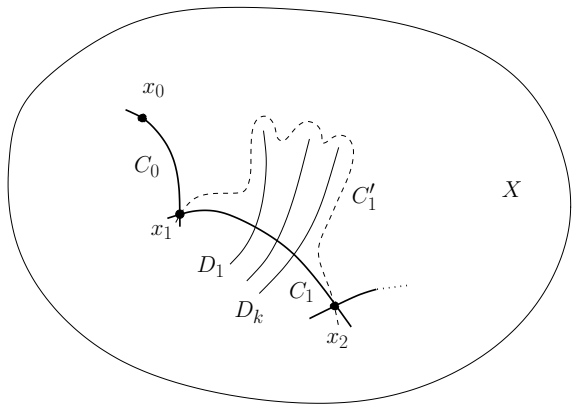}}
\caption{Obtaining a free rational curve through $x_1$ and $x_2$} 
\label{connecting}
\end{figure}

We can deform $C_0'$ to produce many free rational curves 
$D_1,\dots,D_k$ meeting $C_1$. Indeed, view $C_0'$ as the image of a morphism
$f_0:\p^1\to X$, with $f_0(0)=x_1$. Let $U$ be a neighborhood of $[f_0]$ in
$\Hom(\p^1,X)$ parametrizing free morphisms and consider the evaluation 
map $F_U:\p^1\times U\to X$. Since $f_0$ is free, the point
$x_1\in C_0'$ can be moved in any direction, i.e., $F_U(\{0\}\times U)$ contains
a neighborhood of $x_1$ in $X$ (see \cite[II.3.5.3]{Ko96} for a rigorous 
argument). In particular $F_U(\{0\}\times U)$ contains a 
neighborhood of $x_1$ in $C_1$. 
Since $C_0$ meets $X^{\rm free}$, we can choose deformations $D_1,\dots,D_k$ 
that intersect $C_1$ (in distinct points) and also meet $X^{\rm free}$.
By Theorem~\ref{nonfree}, if we take $k$ large enough, we
can deform some subcurve $C_1\cup D_{i_1}\cup \dots \cup D_{i_m}$
of $D\cup D_1\cup \dots \cup D_k$ 
into a rational curve $C_1'$, keeping $x_1$ and $x_2$ fixed.
Since the $D_i$ meet $X^{\rm free}$, 
we can take $C_1'$ to meet $X^{\rm free}$ as well.
(See Figure~\ref{connecting}.)

We do the same with $C_0'$ replaced with $C_1'$, and $C_1$ 
replaced with $C_2$, and so on.



At the end we get   a chain
of free rational curves connecting $x$ and $x'$. 
As in Application \ref{move.free.appl}, we can deform
this chain into a free rational curve, keeping $x$ fixed.
(We really would like to 
deform
this chain into a free rational curve, keeping $x$ and $x'$ both fixed,
but this is not automatic.)

Now consider the evaluation map 
$F:\p^1 \times \Hom(\p^1,X,0\mapsto x)\to X$
sending $(t,[f])$ to $f(t)$.
We have proved that 
 for every $x'\in X^0$ there is an irreducible component
$Z=Z(x')\subset \Hom(\p^1,X,0\mapsto x)$ 
such that $x'$ is in the closure of $F(\p^1 \times Z(x'))$.

The scheme $\Hom(\p^1,X,0\mapsto x)$ has only 
countably many irreducible components, hence
there is a single  irreducible component $Z^*\subset \Hom(\p^1,X,0\mapsto x)$
such that $F|_{\p^1\times Z^*}$ is still dominant.
By Proposition~\ref{smoothpointisfree.vfree}
 we conclude that $f$ is very free for
general $[f]\in Z^*$.
\medskip

Now let us prove (6) $\Rightarrow $ (1$^+$). 

Let $f:\p^1 \to X$ be a very free  morphism.
By Proposition~\ref{smoothpointisfree.vfree}, 
there exists 
a variety $Y$ of dimension $\dim X-1$, a point $y\in Y$ 
and a dominant morphism
$$
F:\PP^1\times Y\to X\qtq{such that}F|_{\PP^1\times \{y\}}=f\qtq{and} 
F|_{\{0\}\times Y}=x_0.
$$

\begin{figure}[hbtp]
\centering

\psfrag{x0}{{\fontsize{28}{20}$x_0$}}
\psfrag{x}{{\fontsize{28}{20}$x_1$}}
\psfrag{C0}{{\fontsize{28}{20}$C_0$}}
\psfrag{C1}{{\fontsize{28}{20}$C_1$}}
\psfrag{C2}{{\fontsize{28}{20}$C_2$}}
\psfrag{Cm}{{\fontsize{28}{20}$C_m$}}
\psfrag{X}{{\fontsize{28}{20}$X$}}

\scalebox{0.35}{\includegraphics{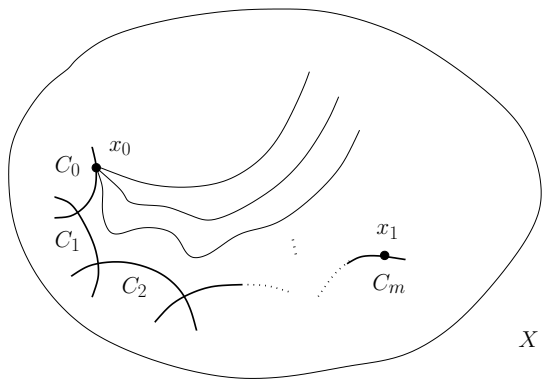}}
\caption{A rational curve moving 
         with one point fixed, and then degenerating into a union 
         of rational curves.}
\label{degeneration}
\end{figure}

If $x_1$ and $x_2$ are in the image of $F:\PP^1\times Y\to X$, we 
have  a length 2 chain of rational curves connecting $x_1$ to $x_0$ and 
then to $x_2$. Otherwise, since $F$ is dominant, 
we still have a  $1$-parameter family $\{f_t\}$ of deformations of $f$
fixing $x_0$ that  degenerates into a tree of rational curves $T_1$
connecting $ x_1$ and $x_0$ (see Figure~\ref{degeneration}). 
Similarly there is a tree of rational curves $T_2$
connecting $ x_2$ and $x_0$.



To see that a $1$-dimensional 
family of rational
curves on $X$ can only degenerate into a union of rational curves, we
represent this family by a morphism $\phi:S\to X$, where
$S$ is a $\p^1$-bundle over a smooth curve $T$. Let $\bar{T}$ be a smooth
compactification of $T$, and let $\bar{S}\to \bar{T}$ be a $\p^1$-bundle 
extending $S\to T$. Let $\bar{\phi}:\bar{S} \map X$ be a rational 
map extending $\phi$. If $\bar{\phi}$ is a morphism, then we are done.
Otherwise we perform a sequence of blow ups of points $\tilde{S}\to \bar{S}$
in order to resolve the indeterminacy of $\bar{\phi}$, and obtain a
morphism $\tilde{\phi}:\tilde{S}\to X$ extending $\phi$. The result 
follows from the fact that blow ups of points only add copies of $\p^1$
to the fibers of $\bar{S}\to \bar{T}$.

\medskip

We complete the proof by showing (6$^+$) $\Rightarrow$ (5).

Assume that for every $x_i$ there is a very free  morphism
 $f_i:\p^1 \to X$  with $f_i(0)=x_i$.
By Proposition~\ref{smoothpointisfree.vfree}, 
there are pointed
 varieties $(Y_i,y_i)$ of dimension $\dim X-1$,  
and  dominant morphisms
$$
F_i:\PP^1\times Y_i\to X\qtq{such that}F_i|_{\PP^1\times \{y_i\}}=f_i\qtq{and} 
F_i|_{\{0\}\times Y_i}=x_i.
$$
We may also assume that $F_i|_{\PP^1\times \{y'\}}$
is very free for every $y'\in Y_i$.

Let $X^*\subset X$ be a dense open set contained in all
the images $F_i(\PP^1\times (Y_i\setminus \{y_i\})$.
Let  $g_0:C_0\cong \p^1\to X$ be a very free morphism intersecting $X^*$
and pick points $c_1,\dots,c_m\in C_0$ such that
$g_0(c_i)\in X^*$. 
By the choice of $X^*$, there are $z_i\in Y_i$
such that  $F_i((\infty, z_i))=g_0(c_i)$.
Set $g_i:=F_i|_{\PP^1\times \{z_i\}}$ (so $g_i(0)=x_i$ and 
$g_i(\infty)=g_0(c_i)$).

Let $C_0\cup C_1\cup \dots \cup C_m$  be a rational comb such that 
$C_0\cap C_i=\{c_i\}$ for $i=1,\dots ,m$. 
Choose identifications $C_i\cong \p^1$ such that $c_i$ is identified
with $\infty$. Define a morphism 
$$
G:C\to X
\qtq{by} G|_{C_i}=g_i, \ 0\leq i \leq m.
$$
Notice that it passes through the points $x_1,\dots ,x_m$.

Like in Application \ref{move.free.appl}, the fact that 
all the  $g_i$ are very free implies that $G:C\to X$ 
can be deformed into a single
very free morphism $\p^1\to X$ passing through all the $x_i$.
 \qed

\medskip

Now that we have established our definition, we can show that
the class of rationally connected varieties enjoys several
good properties.

Let us start with one half of Castelnuovo's theorem.

\begin{thm}\cite{KMM92b} \label{Castelnuovo}
Let $X$ be a smooth, projective, rationally
connected variety over $\c$. Then
$$
H^0(X,\left(\Omega^1_X\right)^{\otimes m})=0\qtq{for all $m\geq 1$.}
$$
\end{thm}

\noindent {\it Proof.} Choose $f:\p^1\to X$ such that $f^*T_X$ is ample
(i.e., $f$ is very free).
Then $f^*\Omega^1_X$ is a direct sum of line bundles with negative degree,
and so is every tensor power. Thus any section of 
$\left(\Omega^1_X\right)^{\otimes m}$
vanishes along $f(\p^1)$. Since $X$ is covered by very free rational 
curves (this is characterization (6$^+$) of Definition-Theorem~\ref{RCV}),
we cannot have a nonzero section.\qed

\medskip

It would be really interesting to know that the converse also
holds, but this is known only in dimension 3 (see \cite[3.2]{KMM92b}).

\begin{thm} \cite{KMM92b} \label{deformation_invariance}
Let $X$ be a variety over $\c$, $S$ a connected $\c$-scheme 
and $p:X\to S$ a smooth projective morphism.
If a fiber $X_s$ is rationally connected for some
$s\in S$, then every fiber of $p$ is rationally connected.
\end{thm}

This is again not hard to prove. We can use
characterization (6) of Definition-Theorem~\ref{RCV} to prove that
being rationally connected is an open condition. Then we use
characterization (1$^+$) to prove that it is also a closed condition.
For a complete proof we refer to \cite[2.4]{KMM92b} or \cite[IV.3.11]{Ko96}.
\medskip

The analog of Noether's theorem was recently established by
Graber, Harris and Starr.

\begin{thm}\cite{GHS01} \label{ghs} 
Let $X$ be a smooth projective variety 
over $\C$ and $f:X\to \PP^1$ a morphism. If the general fiber of $f$ is 
rationally connected, $f$ has a section. \qed \end{thm}

As a consequence of Theorem~\ref{ghs} we have:

\begin{cor} \label{rcfibration}
Let $f:X\to Y$ be a dominant morphism of smooth projective
varieties over $\c$. If $Y$ and the general fiber of $f$ are 
rationally connected, so is $X$.  \end{cor}

\noindent {\it Proof.} 
Let $x_1, x_2 \in X$ be two points such that $f$ is smooth 
with rationally connected fibers over
 $f(x_1)$ and $f(x_2)$.
 We  show that they can be
connected by a chain of rational curves. Since $Y$ is rationally connected,
there is a rational curve joining $f(x_1)$ and $f(x_2)$. 
Represent it by a morphism  $\PP^1\to  Y$ birational onto its image,
let $Z'$ be the unique irreducible component of  $\PP^1\times_YX$
that dominates $\p^1$ 
and let $Z\to Z'$ be a resolution of singularities.
The fibers of $Z'\to \PP^1$ containing $x_1,x_2$ are already smooth, so
we do not need to change them to get $Z$.
Let
$f_Z:Z\to \PP^1$ denote the projection.

By Theorem~\ref{ghs}, $f_Z$ has a section 
$D\subset Z$. Let $C\subset X$ be the projection of $D$ to $X$ and 
$x_1'$ (resp. $x_2'$) be the intersection point of
$C$ with the fiber $X_{f(x_1)}$ (resp. $X_{f(x_2)}$). 
By assumption,  there are  rational curves 
$C_1\subset X_{f(x_1)}$ (resp. $C_2\subset X_{f(x_2)}$) joining $x_1$ 
and $x_1'$ 
(resp. $x_2$ and $x_2'$). Then $C_1\cup C\cup C_2$ is a chain of rational 
curves through $x_1$ and $x_2$ (see Figure~{\ref{figrcfibration}).
\qed
\medskip



\begin{figure}[hbtp]
\centering

\psfrag{X}{{\fontsize{28}{20}$X$}}
\psfrag{Y}{{\fontsize{28}{20}$Y$}}
\psfrag{f}{{\fontsize{28}{20}$f$}}
\psfrag{X1}{{\fontsize{28}{20}$X_{f(x_1)}$}}
\psfrag{X2}{{\fontsize{28}{20}$X_{f(x_2)}$}}
\psfrag{x1}{{\fontsize{28}{20}$x_1$}}
\psfrag{x2}{{\fontsize{28}{20}$x_2$}}
\psfrag{x1'}{{\fontsize{28}{20}$x_1'$}}
\psfrag{x2'}{{\fontsize{28}{20}$x_2'$}}
\psfrag{C}{{\fontsize{28}{20}$C$}}
\psfrag{C1}{{\fontsize{28}{20}$C_1$}}
\psfrag{C2}{{\fontsize{28}{20}$C_2$}}
\psfrag{D}{}
\psfrag{fx1}{{\fontsize{28}{20}$f(x_1)$}}
\psfrag{fx2}{{\fontsize{28}{20}$f(x_2)$}}

\scalebox{0.35}{\includegraphics{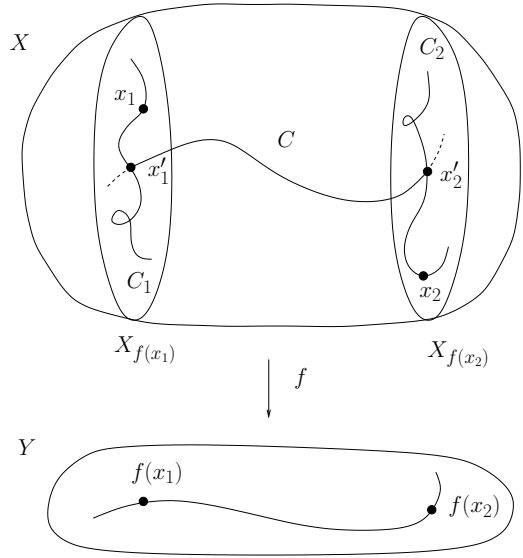}}
\caption{Joining $x_1$ and $x_2$ by a chain of rational curves}
\label{figrcfibration}
\end{figure} 


\begin{say}[Rational connectedness in positive characteristic]

In characteristic $p$, the conditions of Definition-Theorem~\ref{RCV}
are not equivalent, and we have two distinct concepts, rationally
connected varieties, and \emph{\SRC varieties}, which we define below.

Let $X$ be a smooth projective variety defined over an arbitrary uncountable 
algebraically
closed field. The following two conditions are equivalent.
\begin{enumerate}
\item Through two general points of $X$ there passes a rational curve.
\item There exists a variety $Z$ and a morphism $U=Z\times\p^1 \to X$
such that the induced morphism $U\times_Z U \to X\times X$ is 
dominant. (A point $(x_1,x_2)\in X\times X$ is in the image of this morphism
if and only if there is a point $z\in Z$ such that the image of 
$\{z\}\times \p^1$ in $X$ 
contains $x_1$ and $x_2$.)
\end{enumerate}

In positive characteristic these conditions do not imply the
existence of a very free morphism $f:\p^1 \to X$. 

It turns out that the following two conditions are equivalent.
\begin{enumerate} \setcounter{enumi}{2}
\item There is a morphism $f:\p^1\to X$ such that $f^*T_X$ is ample.
\item There exists a variety $Z$ and a morphism $U=Z\times\p^1\to X$
such that the induced morphism $U\times_Z U \to X\times X$ is 
dominant and smooth at the generic point.
\end{enumerate}

This motivates the following definition.
\end{say}

\begin{defn}
Let $X$ be a smooth projective variety over an arbitrary field $k$.
We say 
that $X$ is \emph{\SRC} if there is a morphism (defined over 
$\bar{k}$) $f:\p^1_{\bar{k}}\to X_{\bar{k}}$ such that $f^*T_X$ is
ample.  
\end{defn}

The condition above implies conditions (1) through (6) and also
conditions ($1^+$) and ($6^+$) in Definition-Theorem~\ref{RCV}. 

Over  fields of characteristic $0$, the notions of 
rationally connected and \SRC varieties coincide. However, 
there are examples of varieties in positive characteristic 
satisfying condition~1 above but not 3 (see for instance
\cite[V.5.19]{Ko96}).

We  study \SRC varieties over local fields in section~\ref{nonACF}.


\section{Genus Zero Stable Curves} \label{genus_0_stable_curves}

In the last sections we saw that it is very useful to be able to
deform reducible curves into irreducible ones. 
Frequently, we also need to degenerate 
rational curves so that they become reducible. 
For all of these processes it is very convenient to have a good 
theory that encompasses all genus zero  curves in one object.

One possibility is to work with pairs $(C,g)$,
where $C$ is a connected and possibly reducible curve 
of genus zero (that is, with $h^1(C,\o_C)=0$) and
$g:C\to X$ is a finite morphism.
 Geometrically these objects are
 easy to visualize. A technical problem arises, however, when $C$ has
$m\geq 3$ irreducible components meeting in one point.
Here locally $C$ looks like the $m$ coordinate axes in
$\a^m$. 
There are at least two problems.
\begin{enumerate}
\item The dualizing sheaf $\omega_C$ is not locally free,
making duality less useful.
\item 
The deformation theory of this object is
rather complicated.
\end{enumerate}

\cite{konts} observed that it is  better
to insist that $C$ have only nodes. The price we pay
is that we have to allow $g:C\to X$ to map some
irreducible components of $C$ to points.
This also means that even a reasonably nice curve on
$X$ does not  correspond to a unique pair $(C,g)$.

\begin{exmp} First let $D\subset \p^3$ consist of 3 lines 
$D_1,D_2,D_3$  through
a point $P$ with independent directions. Then $h^1(D,\o_D)=0$. 
The way we represent this object is by a morphism
$g:C\to D$, where $C$ is a rational comb with handle $C_0$ and teeth
$C_1, C_2, C_3$. The morphism $g$ maps $C_i$  
isomorphically onto $D_i$ for $i=1,2,3$ and contracts $C_0$ to $P$. 
We have 3 nodes 
on $C_0$.

Next consider the case when  $D\subset \p^4$ consists of 4 lines 
$D_1,\dots,D_4$  through
a point $P$ with independent directions.  Again $h^1(D,\o_D)=0$. 
As before,  we can represent this object  by a morphism
$g:C\to D$, where $C$ 
is a rational comb with handle $C_0$ and teeth
$C_1, C_2, C_3, C_4$. The morphism $g$ maps $C_i$  
isomorphically onto $D_i$ for $i=1,2,3,4$
and contracts $C_0$ to $P$.
This time we have 4 nodes on $C_0$,
and these four points  have moduli, the cross ratio.

We  also have 3 degenerate cases when 2 of the nodes
on $C_0$ come together. We represent them by 
replacing $C_0$ with two rational curves $C_0'$ and $C_0''$
meeting transversely in a single point. The curves $C_1,C_2,C_3,C_4$ 
are still disjoint from each other. Two of them meet $C_0'$ transversely 
and do not meet $C_0''$, while the others meet $C_0''$
transversely and do not meet $C_0'$. Each $C_i$ maps 
isomorphically onto $D_i$, while $C_0'$ and $C_0''$ are contracted to 
$P$ (see Figure~\ref{combanddegenerations}).
\end{exmp}



\begin{figure}[hbtp]
\centering

\psfrag{C0}{{\fontsize{20}{20}$C_0$}}
\psfrag{Di}{{\fontsize{20}{20}$D_1$}}
\psfrag{Dj}{{\fontsize{20}{20}$D_2$}}
\psfrag{Dk}{{\fontsize{20}{20}$D_3$}}
\psfrag{Dl}{{\fontsize{20}{20}$D_4$}}
\psfrag{C0'}{{\fontsize{20}{20}$C_0'$}}
\psfrag{C0''}{{\fontsize{20}{20}$C_0''$}}
\psfrag{Ci}{{\fontsize{20}{20}$C_i$}}
\psfrag{Cj}{{\fontsize{20}{20}$C_j$}}
\psfrag{Ck}{{\fontsize{20}{20}$C_k$}}
\psfrag{Cl}{{\fontsize{20}{20}$C_l$}}

\scalebox{0.45}{\includegraphics{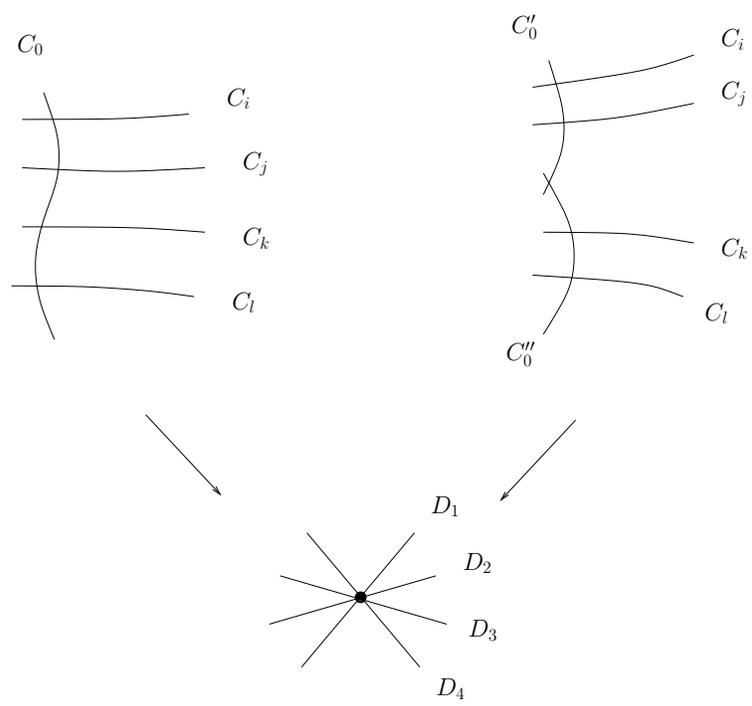}}
\caption{Stable curves representing 4 lines through a point}
\label{combanddegenerations}
\end{figure} 


Next we define stable curves. Then we give some basic properties of the
space parametrizing such objects.

\begin{defn} Let $X$ be a scheme over a field $k$.
A {\it stable  curve over $X$
(defined over $k$)} is a pair $(C,f)$
where 
\begin{enumerate}
\item $C$ is a proper connected curve defined over $k$ having only
nodes,  
\item  $f:C\to X$ is a $k$-morphism, and 
\item $C$ has only finitely many automorphisms that commute with $f$.
Equivalently, the following 2 situations
are excluded:
\begin{enumerate}
\item The curve $C$ has an irreducible component $C'\cong \p^1$ 
that contains at most 2 nodes of $C$ and $f(C')$ is a point.
\item The curve $C$ is a smooth elliptic curve   and 
 $f(C)$ is  a point.
\end{enumerate}
\end{enumerate}
\end{defn}

When $X$ is a point, this is the Deligne-Mumford definition of a 
stable curve \cite{delmum}.

We also need to study cases when
all the curves have 
marked points or they 
 pass through certain points of $X$.
Over an algebraically closed field these are just points,
but in the nonclosed case we should keep track of the
scheme structure on the points.
This is the role of the scheme $P$ in the definition below.

\begin{defn} Let $X$ be a scheme over a field $k$,
 $P$  a smooth zero dimensional $k$-scheme
and $t:P\to X$ a $k$-morphism.
A {\it stable  curve over $X$ (defined over $k$) 
with base  point $t(P)$}
 is a triple $(C,f,\sigma)$
where 
\begin{enumerate}
\item $C$ is a proper connected curve defined over $k$ having only
nodes,  
\item  $\sigma:P\to C$ is a $k$-embedding into the smooth locus of $C$,
\item  $f:C\to X$ is a $k$-morphism such that $t=f\circ \sigma$, and
\item $C$ has only finitely many automorphisms that commute with $f$
and $\sigma$.
\end{enumerate}
\end{defn}

We usually work with families of curves.

\begin{defn}
A {\it family of stable  curves over $X$ with base point $t(P)$}
 is a triple
 $(\cC ,f,\sigma)$
where
\begin{enumerate}
\item
 $\cC \to U$ is a flat, proper family of
 connected curves with only
nodes over a $k$-scheme $U$,
\item  $\sigma:P\times_k U\to \cC$ is an embedding, and
\item     $f:\cC \to X\times_k U$ is a $U$-morphism such that 
$(\cC_u,f_u,\sigma_u)$ is a stable  curve (defined over $k(u)$)
with base point $t_u(P_u)$ for every $u\in U$.
\end{enumerate}
\end{defn}

A basic theorem, due to \cite{konts} 
(in the context of algebraic stacks) 
and \cite{alexeev} 
(in the present setting of projective moduli space), 
asserts that one can 
parametrize all stable curves with a single
scheme. The relevant notion is the concept of
coarse moduli spaces, introduced by Mumford. 
A more general version is proved in section 10.

\begin{thm} \label{moduli}  
Let $X$ be a projective scheme  over a field $k$ with an ample divisor $H$,
$P$  a smooth zero dimensional $k$-scheme 
and $t:P\to X$ a $k$-morphism.
Fix integers $g$ and $d$.
Then there is a projective $k$-scheme $\overline{M}_{g}(X,d,t)$
which is a coarse moduli space for all genus $g$
stable curves $(C,f,\sigma)$ with base point $t(P)$
 such that $f_*[C]\cdot H=d$.\qed
\end{thm}

\begin{say}[What does this mean?] 
We want to say that there is a ``natural''
 correspondence between the set of all stable curves
and the points of $\overline{M}_{g}(X,d,t)$. For our purposes
this means the following two requirements.
\begin{enumerate}
\item There is a bijection of sets
$$ 
\Phi: 
\left\{
\begin{array}{c}
\mbox{genus $g$ stable curves over $\bar k$}\\
 (C,f,\sigma)\\
\mbox{with base point $t(P)$}\\
\mbox{such that $f_*[C]\cdot H=d$}
\end{array}
\right\}
\stackrel{\cong}{\longrightarrow}
\left\{
\mbox{$\bar k$-points of $\overline{M}_{g}(X,d,t)$}
\right\}.
$$
\item 
If $(\cC /U,f,\sigma)$
is a family of degree $d$ stable curves 
with base point $t(P)$, then
there is a unique morphism $M_U:U\to \overline{M}_{g}(X,d,t)$
such that 
$$
M_U(u)=\Phi(\cC_u,f_u,\sigma_u)
\qtq{for every $u\in U(\bar k)$.}
$$
$M_U$  is called the {\it moduli map}.
\end{enumerate}

Given a stable curve $(C,f,\sigma)$ with base point $t(P)$, 
we denote by $[(C,f,\sigma)]$
the corresponding point $\Phi((C,f,\sigma))$ in $\overline{M}_{g}(X,d,t)$.

We usually do not care about the degree, and 
use the scheme $\overline{M}_{g}(X,t)$ which is the 
union of all the 
$\overline{M}_{g}(X,d,t)$. The larger scheme
$\overline{M}_{g}(X,t)$ is not projective since it has
infinitely many connected components.

The requirements above seem somewhat awkward, but one cannot achieve more
in general. For instance, (1) usually fails over
nonclosed fields. The genus zero curves 
$(x_0^2+x_1^2+x_2^2=0)\subset \p^2_{\r}$ and $\p^1_{\r}$ are isomorphic 
over $\c$ but not over $\r$. Over $\q$ there are infinitely many 
nonisomorphic genus zero curves (all isomorphic over $\c$). These
examples ilustrate that $\Phi$ is not  injective over nonclosed
fields. In general, $\Phi$ is also not surjective, but
such examples are harder to get.

\end{say}

We need to know a little about the
infinitesimal description of $\overline{M}_{g}(X,t)$
near a stable curve $(C,f,\sigma)$.
As usual, if $X$ is smooth, this is governed by the
cohomology groups of the twisted pull back of the tangent bundle,
$f^*T_X(-\sigma(P))$. 
The moduli space is usually very complicated,
but its local structure is nice
in some cases.

Trouble comes from two sources.
First, if $H^1(C,f^*T_X(-\sigma(P))$ is large then
one can say little beyond
a dimension estimate, as in Theorem~\ref{ksm}.3.
Second, if 
$C$ has automorphisms commuting with $f$ and $\sigma$,
then $\overline{M}_{g}(X,t)$ looks like the quotient of a
smooth variety by a finite group.
This can be quite complicated, especially
in positive characteristic. 
This accounts for the somewhat messy formulation
of Theorem~\ref{Mg.near.free.thm}.5.
Fortunately, this seemingly cumbersome form is
more useful for many purposes. 
A proof is given in section 10.

\begin{thm} \label{Mg.near.free.thm}
Let $X$ be a smooth projective variety over a field $k$.
Let $P$ be a smooth zero dimensional $k$-scheme  and 
$t:P\to X$  a $k$-morphism. Set $|P|:=\dim_k\o_P$.
Let $(C,f,\sigma)$ be a    stable curve 
of genus $g$  with base point $t(P)$. 
Assume $H^1(C,f^*T_X(-\sigma(P)))=0$.
Then 
\begin{enumerate}
\item
There is a unique irreducible component (defined over $k$) 
$\comp(C,f,\sigma)\subset \overline{M}_{g}(X,t)$ passing through 
$[(C,f,\sigma)]$.
\item 
$\comp(C,f,\sigma)$ is geometrically irreducible.
\item
There is a dense open subset of $\comp(C,f,\sigma)$ that corresponds to
irreducible   curves.
\item 
$\comp(C,f,\sigma)$ has dimension 
$$
-(f_*[C]\cdot K_X)+(\dim X-3)(1-g)-|P|(\dim X-1).
$$
\item 
There is a smooth $k$-scheme $U$ with a rational point $u\in U(k)$ 
and a family of stable curves 
$$
(\cC/U, F_U:\cC\to X\times U, \sigma_U:P\times U\to \cC)
$$ 
such that $(\cC_u,F_u,\sigma_u)\cong (C,f,\sigma)$ and
the corresponding moduli map
$$
M_U:(U,u)\to (\comp(C,f,\sigma),[(C,f,\sigma)])
$$
is  quasifinite and  open.
\end{enumerate}
\end{thm}

In fact, if  $C$ has no automorphisms commuting with $f$
and $\sigma$, we can choose $M_U$ to
be an open embedding.


\section{Finding Rational Curves over Nonclosed Fields}
\label{nonACF}

Let $X$ be a smooth, projective, \SRC variety over a field $k$.
By definition, $X$ has plenty of rational curves over $\bar k$,
but we know very little about rational curves over $k$.

It is possible that $X(k)=\varnothing$, and then
we certainly do not have any morphism  $\p^1\to X$
defined over $k$.

Thus assume there exists a point $x\in X(k)$. Can we find a
rational curve $f:\p^1\to X$ defined over $k$ and passing through $x$?
It is quite unlikely that such a curve always exists,
but no counterexample is known.
Here we discuss a weaker version of the problem,
which allows us to settle the above question for
certain fields.

Consider the space $\overline{M}_0(X,\spec k\to x)$ of genus zero stable 
curves through $x$
and let $Z_i$ be its irreducible components over $\bar k$.
The Galois group $\gal(\bar k/k)$ acts on the set of all
components. If there is a free rational curve (defined over $\bar k$)
$f:\p^1_{\bar k}\to X_{\bar k}$, $f(0)=x$,
then  $[(\p^1,f,\spec k\to 0)]$ lies in a unique $Z_i=\comp(f)$ by 
Theorem~\ref{Mg.near.free.thm}. If, in addition, $f$ is 
defined over $k$, then  $\comp(f)$ is a geometrically
irreducible component of $\overline{M}_0(X,\spec k\to x)$ 
defined over $k$.

In order to answer the above question, we do the following.
First we construct geometrically
irreducible components $Z_i$ of $\overline{M}_0(X,\spec k\to x)$ that are
defined over $k$. This is done over any field. 
For local fields
we can then find $k$-points in $Z_i$ corresponding to irreducible 
curves through $x$. For arbitrary fields
this seems rather difficult.

The main method is to assemble combs from given morphisms.
This is what we define next.

\begin{defn} Let $X$ be a $k$-variety and $x\in X(k)$.
Assume that we have $\bar k$-morphisms (not necessarily defined over $k$)
$f_i:\p^1_{\bar k}\to X_{\bar k}$ such that $f_i(0)=x$ for $i=1,\dots,n$.

Let $C:=C_0\cup \dots \cup C_n$ be a rational comb 
(over $\bar k$) with handle $C_0$.
Choose identifications $C_i\cong \p^1$ such that the node
is identified with $0$. Define a morphism
$$
F:C\to X
\qtq{by} F|_{C_i}=f_i\qtq{and}F(C_0)=\{x\}.
$$ 
The pair $(C,F)$ is called a {\it comb assembled from the $f_i$}.
If $n\geq 3$, $(C,F)$ is a stable
curve over $X_{\bar k}$, provided none of the $f_i$
are constant.

A comb assembled from the $f_i$ is not unique,
but all the combs assembled from the $f_i$
form an irreducible family
(parametrized by the space of $n$ points on $C_0\cong \p^1$).

It is natural to think of the $f_i$
as $\bar k$-points of $\overline{M}_0(X,\spec k\to x)$.
In this context, we pick an extra point $p_0\in C_0(\bar k)$
distinct from the nodes of $C$ and
set $\sigma:\spec \bar k\to p_0$. This makes
$(C,F,\sigma)$ into a 
stable curve (defined over $\bar k$) with base point $x$,
provided that $n\geq 2$ and none of the $f_i$
are constant.
We call this triple a
{\it base pointed comb assembled from the $f_i$}.

Assume next that the
$$
f_i:\p^1\to X,\ f_i(0)=x,\qtq{$i=1,\dots, n$,}
$$
form a conjugation invariant set over $k$. 
We would like to assemble all the $f_i$ into an
element of $\overline{M}_0(X,\spec k\to x)$ that is defined over $k$.
The construction above has to be carried out carefully: even if we choose
$C$ to be a comb over $k$, the morphism $F$ is not defined over $k$ in
general. In order to achieve this, we need that the induced scheme structure 
on the subset of nodes $\{P_i=C_0\cap C_i\} \subset C_0$ matches the one on
$\{[(C_i,f_i)]_{i=1,\dots,n}\}\subset \overline{M}_0(X,\spec k\to x)$. 
This is the role of
the scheme $T$ below. 

View each $f_i$ as a $\bar k$-point of $\overline{M}_0(X,\spec k\to x)$.
This identifies the collection $\{f_i\}$ 
with a  zero dimensional smooth $k$-subscheme 
$T\subset \overline{M}_0(X,\spec k\to x)$.
(See Remark~\ref{insep.rem.comb} for some further clarification
in positive characteristic.)
The $f_i$ correspond to a $k$-morphism $f_T:\p^1_T\to X$
which maps $\{0\}\times T$ to $x$. 
This construction is interesting only when none of the $f_i$ 
corresponds to a $k$-point of $\overline{M}_0(X,\spec k\to x)$, 
so we assume that $T$ has no $k$-points.

Let  $C_0= \p^1$ and fix an embedding $j:T\into C_0$.
This determines a rational comb $C$ over $k$ with teeth $\p^1_T$
whose nodes coincide with $j(T)$.
Now define a morphism
$$
F:C\to X
\qtq{by} F|_{\p^1_T}=\{f_i\}\qtq{and}F(C_0)=\{x\}.
$$ 
Notice that $C$ and $F$ are both defined over $k$. Since $T$ has no 
$k$-points,
we can pick a point $p_0\in C_0(k)$ distinct from the nodes of $C$
and set $\sigma:\spec k \to p_0$.
The triple
$(C,F,\sigma)$ is called a {\it pointed $k$-comb assembled from the $f_i$}.
If $n\geq 2$, $(C,F,\sigma)$ is a stable
curve over $X$ with base point $x$, provided none of the $f_i$
are constant.

As before, a $k$-comb assembled from the $f_i$ is not unique,
but all $\bar k$-combs assembled from the $f_i$ 
form an irreducible family
defined over $k$
(parametrized by the space of all $\bar k$-embeddings  
$T\cup \spec k \into \p^1$).
\end{defn}

\begin{prop} 
\label{comb.free.prop}
Let $k'/k$ be a finite, separable field extension.
 Let $X$ be a smooth projective variety over $k$ and let $x\in X(k)$
be a point. Let 
$$
f_i:\p^1\to X,\ f_i(0)=x,\qtq{$i=1,\dots, n$,}
$$
be a set of nonconstant free $k'$-morphisms  together defined  over $k$.
Let $(C,F,\sigma)$ be a pointed $k$-comb assembled from the $f_i$,
so that $(C,F,\sigma)$ is a stable curve defined over $k$ with base point $x$.

Then there is a smooth $k$-scheme $U$ with a $k$-point $u\in U(k)$ 
and a family of stable curves 
$$
(\cC/U, F_U:\cC\to X\times U, \sigma_U:P\times U\to \cC)
$$ 
such that $(\cC_u,F_u,\sigma_u)\cong (C,F,\sigma)$ and
the corresponding moduli map
$$
M_U:(U,u)\to (\comp(C,F,\sigma),[(C,F,\sigma)])
$$
is  quasifinite and  open.


Moreover, the irreducible component 
$\comp(C,F,\sigma)\subset \overline{M}_0(X,\spec k\to x)$
is geometrically irreducible and independent of the choice of the comb.
\end{prop}

\noindent {\it Proof.} Denote by $C_0$ the handle of $C$ and by 
$C_1,\dots, C_n$ its teeth. Let $p_0=\sigma(\spec k)\in C_0$.

Set $E=F^*T_X(-p_0)$. Since $F$ maps $C_0$ to
a point, $f^*T_X|_{C_0}=\o_{C_0}^{\oplus n}$. Thus
$H^1(C_0,E|_{C_0})=0$. For any other $C_i$,
$E|_{C_i}(-1)\cong f_i^*T_X(-1)$, and this has vanishing $H^1$
since $f_i$ is free. Thus
$H^1(C,F^*T_X(-p_0))=0$
by Exercise~\ref{h1=0.exrc}. 
The first assertion then follows from Theorem~\ref{Mg.near.free.thm},
and so does geometric irreducibility of $\comp(C,F,\sigma)$.

Theorem~\ref{Mg.near.free.thm} also implies that
each $[(C,F,\sigma)]$ lies in a unique
irreducible component  
$\comp(C,F,\sigma)\subset \overline{M}_0(X,\spec k\to x)$.
All $\bar k$-combs assembled from the $f_i$ 
form an irreducible family defined over $k$.
Hence they all lie in the same irreducible component of 
$\overline{M}_0(X,\spec \bar k\to x)$, and
$\comp(C,F,\sigma)$  does not depend on the choice of the
comb. \qed

\begin{rem}\label{insep.rem.comb}
In positive characteristic we have to deal with the possibility that
$f_i:\p^1\to X$ is not defined over a separable field extension.
Since $f_i$ is free, its irreducible component
in $\Hom (\p^1,X, 0\mapsto x)$ is generically smooth,
thus it contains a dense set of points defined over
separable field extensions.
\end{rem}

As an application,  we
prove a  result about \SRC varieties over  local fields.

\begin{defn}
For our purposes  {\it local fields} are
\begin{enumerate}
\item Finite degree extensions of the $p$-adic fields $\q_p$,
\item Power series fields over finite fields, $\f_q((t))$, 
\item $\r$ and $\c$.
\end{enumerate}
For all of these fields we have an analytic implicit function
theorem, thus if $X$ is a variety of dimension $n$ over a local field $k$,
and $x\in X(k)$ is a smooth point, then  a neighborhood of
$x$ in the analytic variety $X(k)$ is biholomorphic
to a neighborhood of $0$ in $k^n$.

In particular, if $X$ has {\em one} smooth $k$-point then
$k$-points are Zariski dense.
\end{defn}

\begin{thm} \cite{kol99}
Let $k$ be a local field
and $X$ a  smooth, projective, \SRC variety defined over $k$.
Let $x\in X(k)$ be a point.

Then there is a nonconstant free morphism defined over $k$
$$
f:\p^1\to X\qtq{such that} f(0)=x.
$$
\end{thm}

\noindent {\it Proof.} By assumption and by Remark~\ref{insep.rem.comb},
 there is such a morphism $g$ over a
separable extension of finite degree $k'/k$. Use 
$g$ and all of its conjugates to assemble a
$k$-comb $(C,F, \spec k\to p_0)$ that is a stable curve
with base point $x$.

By Proposition~\ref{comb.free.prop}, there is a smooth
pointed $k$-scheme $(U,u)$ and a dominant, quasifinite
morphism 
$$
(U,u)\to (\comp(C,F, \spec k\to p_0),[(C,F,\spec k\to p_0)]).
$$

Since $k$ is a local field, $U(k)$ is Zariski dense
in $U$, thus $k$-points are
Zariski dense in $\comp(C,F, \spec k\to p_0)$.
By Theorem~\ref{Mg.near.free.thm}.3, 
there exists a $k$-point $[(C',F',\sigma')]$ that corresponds to
an irreducible curve $F':C'\to X$ through $x$. Since
$C'$ comes with a $k$-point $\sigma'(\spec k)$, 
there is an isomorphism $\p^1\to C'$ defined over $k$.
So we obtain a nonconstant $k$-morphism $\p^1\to X$
whose image passes through $x$.\qed

\medskip

The construction of Proposition~\ref{comb.free.prop} can be made
simultaneously for all points. This is the content of the next 
proposition.

\begin{prop}\label{good.maps.thru.pts} 
Let $X$ be a smooth, projective, \SRC variety over a field $k$.

There is an $X$-scheme $U$ 
and a family of stable curves 
$$
(\cC/U, F_U:\cC\to X\times U, \sigma_U:\spec k \times U\to \cC)
$$ 
such that for each $x\in X$ 
the corresponding moduli map
$$
M_{U_x}:U_x\to \overline{M}_0(X_{k(x)},\spec k(x)\to x)
$$
is  quasifinite and open and maps $U_x$ onto an open subset of
$\overline{M}_0(X_{k(x)},\spec k(x)\to x)$
parametrizing free stable curves. 
\end{prop}

\noindent {\it Proof.}
Let $x_g\in X$ be the generic point. As in Proposition~\ref{comb.free.prop},
we obtain a $k(X)$-scheme $U$ 
and a family of stable curves 
$$
(\cC/U, F_U:\cC\to X_{x_g}\times U, \sigma_U:\spec k(X) \times U\to \cC)
$$ 
such that for each $x\in X$ 
the corresponding moduli map
is  quasifinite and open and maps $U$ onto an open subset of
$\overline{M}_0(X_{k(X)},\spec k(X)\to x_g)$
parametrizing free stable curves.

By specialization, for $x$ in an open set $X^0\subset X$
we obtain a family of stable curves 
$$
(\cC/U_x, F_{U_x}:\cC_{U_x}\to X_{k(x)}\times U_x, \spec k(x) \times U_x\to \cC_{U_x})
$$ 
such that the corresponding moduli map
is  quasifinite and open and maps $U_x$ onto an open subset of
$\overline{M}_0(X_{k(x)},\spec k(x)\to x)$
parametrizing free stable curves.

Repeat the argument with $x_g$ replaced with the generic points
of $X\setminus X^0$.
Continuing in this way we eventually  take care of all
points of $X$.

It is quite likely that one can choose $U$ to be irreducible
and $U\to X$ smooth but the above construction
gives neither.\qed


\section{Moduli spaces of stable curves}

The aim of this rather technical section is to 
prove 
 Theorems \ref{modexists.general} and 
\ref{modexists.prop.general},
which are generalizations of
Theorems \ref{moduli} and \ref{Mg.near.free.thm}.
The main outlines of the proofs follow the general philosophy
of constructing moduli spaces. The closest example is
the special case explained in \cite{FuPa}. There are
two key differences in the setup.

First, in \cite{FuPa} the marking of curves over $\c$ was done by
choosing $n$ smooth points of a curve.
In the general setting it is more natural to let the marking
be an embedding of a zero dimensional scheme $Q$
into the smooth locus of a curve.
The theory
is interesting only when $Q$  has embedding dimension at most $1$,
and the global theory works best only if it is smooth.

Second, we also want to account for the possibility that some
of the marked points have a fixed image. For existence
questions this is not much different, but  the local structure
has to be worked out by a somewhat different construction.

We also drop the ``convexity'' condition of \cite{FuPa},
and replace it with the assumption $H^1(C,f^*T_{X_k}(-\sigma(P_k))=0$,
which is essentially the local version of it.

First we generalize the notion of stable curves, defined in 
Section~\ref{genus_0_stable_curves},
to the relative setting. This means that instead of
working with schemes over a field $k$, we fix a Noetherian scheme $S$,
and work with schemes over $S$.

\begin{defn}\label{stable.gen.defn}
Fix a Noetherian scheme $S$.
Let $X$ be a scheme over $S$, 
let $P\to S$ and $Q\to S$ be  finite, flat  morphisms
and let $t:P\to X$ be a morphism.
A {\it stable  curve over $X$ with base point $t(P)$ and marking $Q$}
 is a quadruple $(C,f,\sigma,\rho)$
where
\begin{enumerate}
\item
 $C\to U$ is a flat, proper family of
 connected curves with only
nodes over an $S$-scheme $U$,
\item   $f:C\to X\times_S U$ is a $U$-morphism,
\item  $\sigma:P\times_S U\to C$  and $\rho:Q\times_S U\to C$ are embeddings
into the smooth locus of $C\to U$
with disjoint images,
\item $f\circ\sigma$ is the morphism $P\times_S U\to X\times_S U$
induced by $t$, and
\item    for every point $\spec k\to U$, $C_k$ has only
finitely many automorphisms commuting with 
$f_k$ and fixing the closed points
of  $\sigma_k(P_k)$ and $\rho_k(Q_k)$.
\end{enumerate}

Note: Condition 5 is equivalent to saying that every rational
component of $C_{\bar k}$ 
that is mapped to a point by $f_{\bar k}$
has at least 3 nodes or marked points
on it. We allow marking by a nonreduced $P_k$,
but a node plus a single point in $P_k$ with a length 3 nilpotent structure 
is not enough for stability, although
it is enough to kill all commuting automorphisms.
 The reason for imposing this stronger condition  
is that without it we get a nonseparated moduli problem.

As a simple example, consider a smooth curve $B$
of genus $g\geq 2$ defined over a field $k$.
Set $X=B$ and $P=\spec k[t]/(t^3)$.
Take  $U=\a^1$ and consider the family $C=B\times U\to U$
with a section
$S=\{b\}\times U$. Denote by $\cI_S$ the ideal sheaf of $S$ in $C$. 
Let $\sigma$ map $P\times U$ onto the
subscheme of $C$  with ideal sheaf $\cI_S^3$. 
Now let $f':C'\to X\times U$ be the blow up of 
$(b,0)$, let $\sigma ':P\times U \to C'$ be the map induced by $\sigma$
and consider the triple $(C',f',\sigma ')$.  
This shows that the identity map $B\to X$
with $P\into B$ degenerates to  $B'\to X$ 
 where $B'=B\cup \p^1$ and $P$ is embedded into $\p^1$.

\end{defn}

In the sequel we assume that $\deg(P/S)$ and $\deg(Q/S)$ are
actually constants. This is the case if $S$ is connected,
so this is not a significant restriction.
We use the shorthand $|P|=\deg(P/S)$ and $|Q|=\deg(Q/S)$.
\medskip

The next theorem is a generalization of Theorem~\ref{moduli}.

\begin{thm} \label{modexists.general}
Fix a Noetherian base scheme $S$. 
Let $X$ be a projective scheme over $S$ with an ample divisor $H$.
Let $P\to S$  and $Q\to S$ be finite and flat morphisms over $S$
of relative embedding dimension at most $1$ 
and let $t:P\to X$ be a morphism.
Fix integers $g$ and $d$.
Then there is a separated algebraic space
${\overline M}_{g,Q}(X,d,t)$
of finite type over $S$
which is a coarse moduli space for all genus $g$
stable curves $(C,f,\sigma,\rho)$
with base point $t(P)$ and marking $Q$ 
 such that $\deg_Cf^*H = d$.

Moreover, ${\overline M}_{g,Q}(X,d,t)\to S$ is projective if $P\to S$  
and $Q\to S$ 
are smooth over $S$.
\end{thm}

It is quite likely that ${\overline M}_{g,Q}(X,d,t)\to S$ is
always quasi-projective, but the proof does not give this.
The reason is that the projectivity criterion
of \cite{koll-proj} used in \cite{FuPa}
does not have a good quasi-projective analog.

As explained in \cite[5.1]{FuPa}, one can 
use $mH$ to embed
$X$ into $\p^M_S$ and deduce the existence
of ${\overline M}_{g,Q}(X,d,t)$ from  the existence
of ${\overline M}_{g,Q}(\p^M_S,md,t)$.
Thus we can restrict ourselves to the 
case when $X$ is smooth over $S$. This is the only case in which
the local structure of ${\overline M}_{g,Q}(X,d,t)$
is easy to describe.
Such a description is given in the next theorem, 
which is a generalization of
Theorem~\ref{Mg.near.free.thm}.
As before, we work with the space ${\overline M}_{g,Q}(X,t)$,
which is the union of all the ${\overline M}_{g,Q}(X,d,t)$, and 
we denote by $[(C,f,\sigma,\rho)]$ the point in 
${\overline M}_{g,Q}(X,t)$ corresponding to a stable curve 
$(C,f,\sigma, \rho)$.

\begin{thm} \label{modexists.prop.general}
Let $S$ be a normal  scheme and
 $X$  a smooth projective scheme over $S$.
Let $P\to S$ and $Q\to S$  be finite and flat  $S$-schemes 
of relative embedding dimension at most $1$ and let
$t:P\to X$ be a morphism. Let $k$ be a field,
$\spec k\to s\in  S$ a point  and 
 $$
(C\to \spec k,f:C\to X_k, \sigma:P_k\to C, \rho:Q_k\to C)
$$ 
 a   stable curve 
of genus $g$ and degree $d$ with base point $t_k(P_k)$
and marking $Q_k$ over $k$. 
 Assume $H^1(C,f^*T_{X_k}(-\sigma(P_k))=0$. Then 
\begin{enumerate}
\item
There is a unique irreducible component 
$\comp(C,f,\sigma,\rho)\subset {\overline M}_{g,Q}(X,t)$
 passing through $[(C,f,\sigma,\rho)]$.
\item If $k$ is the residue field of  $s$ in $S$, then the unique 
irreducible component of the fiber of ${\overline M}_{g,Q}(X,d,t)\to S$
over $s$ containing $[(C,f,\sigma,\rho)]$ is 
geometrically irreducible.
\item
There is a dense open subset of $\comp(C,f,\sigma,\rho)$ that corresponds to
irreducible
 stable curves.
\item $\comp(C,f,\sigma,\rho)$ has dimension 
$$
-(f_*[C]\cdot K_X)+(\dim X-3)(1-g)+\dim S +|Q|-|P|(\dim X-1).
$$
 \item Assume $k/k(s)$ is a finite field extension,
where $k(s)$ denotes the residue field of $s$ in $S$. Then
there is an $S$-scheme $U$, smooth over $S$, with a point
$\spec k \stackrel{\cong}{\to}  u\in U$ 
and a family of stable curves 
$$
(C/U, F_U:C\to X_U, \sigma_U:P_U\to C,
 \rho_U:Q_U\to C)
$$ 
such that 
$(C_u,F_u,\sigma_u,\rho_u)\cong (C,f,\sigma,\rho)$
and the corresponding moduli map
$$
M_U:(U,u)\to (\comp(C,f,\sigma,\rho), [(C,f,\sigma,\rho)])
$$
is  quasifinite and  open.
\end{enumerate}
\end{thm}

If  $C$ has no automorphism commuting with $f$, $\sigma$
and $\rho$, we can choose $M_U$ to
be an open embedding.

\begin{say}[Outline of the proof]
\label{moduli.reduce.gen.rems}

Let $(C,f,\sigma,\rho)$ be a stable curve  with base point $t(P)$,
where $h:C\to U$ is a $U$-scheme and $U$ is a connected $S$-scheme.
The starting point, as in \cite[2.3]{FuPa},
is the observation that $\omega_{C/U}(\sigma(P)+\rho(Q))\otimes f^*\o_X(3H)$
is $h$-ample on $C/U$.  Let $L$ be the above line bundle
raised to the  $(2+|P|)$-th power. Then
$L(-\sigma(P))$ is $h$-very ample  and 
$R^1h_*L=R^1h_*L(-\sigma(P))=0$. Let $N+1$ be
the rank of the vector bundle $h_*L$
(where $N$ depends on the genus of $C$, $|P|$, $|Q|$ and $d$).
After replacing $U$ with an open affine cover, we may assume that
$h_*L\cong \ou_U^{N+1}$. Thus we have an embedding
$\delta:C/U\into \p^N_U$. This is of course not unique,
 but defined only up to automorphisms of $\p^N_U$.

Let us now look at the diagonal map
$$
(f,\delta):C/U\into X\times_U\p^N.
$$
Our strategy has 3 main parts:
\begin{enumerate}
\item  Identify all possible lifts $(f,\delta)$ with some subset of a 
Hilbert scheme. This is straightforward if $P=\varnothing =Q$,
but needs an extra twist otherwise.
\item Construct $ {\overline M}_{g,Q}(X,d,t)$ as the quotient of the above
subset of a
Hilbert scheme by $PGL(N+1)$.
\item Read off properties of $ {\overline M}_{g,Q}(X,d,t)$
through the properties of the Hilbert scheme.
\end{enumerate}

\end{say}

\begin{say}[Construction of a related Hilbert scheme]{\ }
\label{J.const.say}

In order to deal with the base point $t(P)$, let us fix
an embedding $t':P_U\into \p^N_U$
that is nondegenerate
(that is, over any algebraically closed field $K$,
the image $t'(P_K)$ spans a subspace of dimension $|P|-1$).
First we consider only those liftings $\delta$ such that
$t'=\delta\circ\sigma$. In this case the image
of $(f,\delta)$ passes through the image of $(t,t')$,
so it is more natural to view it as an embedded curve
in the scheme $B_{t'}=B_{t'}( X\times_U\p^N)$
obtained by blowing up the image of $(t,t')$.

We need to check a few things since
the scheme $B_{t'}( X\times_U\p^N)$ is singular if $P$ is not smooth.
Let us consider a single fiber of $ X\times_U\p^N\to U$.
We assume that  $P$ has embedding dimension at most $1$ 
and $(t,t')$ is an embedding, 
thus in suitable local analytic coordinates 
$z_1,\dots,z_k$  (where $k=\dim X/U+N$)
we can writte
$$
C=(z_1=\cdots=z_{k-1}=0)\qtq{and} P=(z_1=\cdots=z_{k-1}=z_k^d=0).
$$
The key chart of the blow up is given by coordinates
$z'_1,\dots,z'_{k-1},z_1,\dots,z_k$ and relations
$$
z_i=z'_iz_k^d\qtq{for $i=1,\dots,k-1$.}
$$
This chart is smooth with coordinates $z'_1,\dots,z'_{k-1},z_k$
and the birational transform of $C$ is
$C'=(z'_1=\cdots=z'_{k-1}=0)$. Thus 
$B_{t'}$ is smooth along $C'$.
(All the other local charts of the blow up are singular.)

Inside $\hilb(B_{t'}( X\times_U\p^N))$ consider the
open subset  $W=W_{U,t',g}(X,d,t)$
parametrizing 1-dimensional subschemes $D$ such that
\begin{enumerate}
\item $D$ is a connected curve of genus $g$ 
contained in the smooth locus of $B_{t'}$ and with at worst nodes,
\item the natural projection $D\to \p^N$ is a nondegenerate embedding,
\item  $D$ is smooth at the points of $E$, and 
$E\cap D$ projects isomorphically onto the image of $(t,t')$. 
\end{enumerate}
We still have to add in the marking by $Q$.
To do this consider the universal family
${\mathbf D}\to W_{U,t',g}(X,d,t)$ and
let $Q_W$ denote the pull back of $Q$ to $W_{U,t',g}(X,d,t)$.
Consider the open subset 
$$
J_{U,t',g,Q}(X,d,t)\subset \Hom_{W}(Q_W,{\mathbf D}),
$$
parametrizing curves $D$ and homomorphisms $\rho:Q\to D$
such that
\begin{enumerate}\setcounter{enumi}{3}
\item $\rho$ is an embedding into the smooth locus and its
 image is disjoint from $E\cap D$,
\item the natural projection $D\to X$ is a
stable curve of degree $d$ with base point $t(P)$ and marking $Q$,
\item  on every irreducible component of every geometric
fiber of $D\to U$,
the degree of  $\o_{\p^N}(2+|P|)$ is the same as  the 
degree of  
$$
\omega_{D/U}\otimes \o_X(3H)\otimes \o_{B_{t'}}(E)\otimes \o_D(\rho(Q)).
$$
(The various pull backs to $D$ have been omitted
from the notation.)
\end{enumerate}

Since we are looking only at those
morphisms in $\Hom_{W}(Q_W,{\mathbf D})$
that land in the smooth locus of ${\mathbf D}\to W$,
we conclude from the general case of Theorem~\ref{ksm}
given in \cite[I.2.17]{Ko96}
that
$$
J_{U,t',g,Q}(X,d,t)\to \hilb(B_{t'}( X\times_U\p^N))
$$
is smooth of relative dimension $|Q|$.

If the genus of $C$ is zero, this is the
subscheme we want. In general, based on our construction we also
want to achieve that 
\begin{enumerate}\setcounter{enumi}{6}
\item  the pull back of  $\o_{\p^N}(2+|P|)$ to $D$
is isomorphic to
  $$
\omega_{D/U}\otimes\o_X(3H)\otimes \o_{B_{t'}}(E)\otimes \o_D(\rho(Q)).
$$
\end{enumerate}
As explained in \cite[2.2]{FuPa}, the latter is a locally closed
condition and defines $J^*_{U,t',g,Q}(X,d,t)\subset J_{U,t',g,Q}(X,d,t)$,
which is open if $g=0$ and has codimension $g$ 
in general.

Ultimately we do not want to fix the embedding $t':P\to \p^N$.
To achieve this, let $\Hom^1(P,\p^N)$
denote the open subset of all nondegenerate embeddings.
We have a universal embedding
$$
T: P \times_S\Hom^1(P,\p^N)\to  \p^N\times_S \Hom^1(P,\p^N).
$$
Out of $T$ and $t$ we get a diagonal embedding
$$
(t,T):P \times_S\Hom^1(P,\p^N)\to X\times_S \p^N\times_S \Hom^1(P,\p^N).
$$
Observe that the automorphism group $\aut(\p^N)$ acts
on both sides and $(t,T)$ is equivariant.
Let $B\to X\times_S \p^N\times_S \Hom^1(P,\p^N)$ denote the
blow up of the image.
We can now construct
$$
W_{g}(X,d,t):=W_{\Hom^1(P,\p^N),(t,T),g}(X,d,t)
\subset  \hilb( X\times_S \p^N\times_S \Hom^1(P,\p^N))
$$
and 
$$
J_{g,Q}(X,d,t):=J_{\Hom^1(P,\p^N),(t,T),g,Q}(X,d,t)\to
W_{g}(X,d,t).
$$
The pull back of $J_{g,Q}(X,d,t)$ by the morphism $U\to S$ 
and the embedding $t'$ is
 precisely our earlier $J_{U,t',g,Q}(X,d,t)$.
As before we get that:

\begin{claim}\label{Q.gives.smooth.over}
 With the above notation,
$$
J_{g,Q}(X,d,t)\to W_{g}(X,d,t)
$$
is smooth of relative dimension $|Q|$.\qed
\end{claim}

The automorphism group $\aut(\p^N)$ acts on $B$,
hence it also acts on $\hilb(B)$ and
on $J_{g,Q}(X,d,t)$. From our construction
we have the key observation:

\begin{claim}\label{M=orbits.claim} Let $K$ be any algebraically closed field
and $\spec K\to S$ a morphism. 
\begin{enumerate}
\item If $g=0$,
there is a one-to-one correspondence
between
 the $PGL(N+1,K)$ orbits on
$J_{g,Q}(X_K,d,t_K)$
and 
genus zero stable curves over $K$ of degree $d$ with base point $t(P_K)$ and
marking $Q_K$.
\item If $g\geq 0$,
there is a one-to-one correspondence
between
 the $PGL(N+1,K)$ orbits on
$J_{g,Q}(X_K,d,t_K)$
and 
genus zero stable curves over $K$ of degree $d$ with base point $t(P_K)$ and
marking $Q_K$
plus the choice of a line bundle which has degree zero on every
irreducible component of the curve.
\qed
\end{enumerate}
\end{claim}

\begin{rem} The slight mystery is why we need $K$ algebraically closed.
This is already visible in the much simpler example
of the hyperelliptic plane curves $ay^2=f(x)$,
where we can even  fix $f$ for simplicity.
The multiplicative group $K^*$ acts on the family by $(x,y)\mapsto (x,ty)$
and over an algebraically closed field we have a single orbit.
On the other hand, if $K=\q$, the $\q^*$-action on the $\q$-points has infinitely
many orbits, parametrized by the square free integers.
Upon deeper analysis one sees that 
all of this is caused by the automorphism $(x,y)\mapsto (x,-y)$
of the hyperelliptic plane curves.
\end{rem}
\end{say}

This completes the first step of the plan, and
next we need to take the quotient $J(X,d,t)/\aut(\p^N)$.

\begin{say}[Construction of ${\overline M}_g(X,d,t)$ as a quotient space]{\ }

As we noted above, we would like to find a morphism
$$
q:J_{g,Q}(X,d,t)\to {\overline M}_{g,Q}(X,d,t)
$$
with the property that for
any algebraically closed field $K$ and for any morphism
 $\spec K\to S$  the fibers of
$$
q_K:J_{g,Q}(X_K,d,t_K)\to {\overline M}_{g,Q}(X_K,d,t_K)
$$
are precisely the $PGL(N+1,K)$ orbits on $J_{g,Q}(X_K,d,t_K)$.

 At this point it is very convenient to use the
general quotient theorems about algebraic group   actions
established in \cite{koll-quo} and \cite{ke-mo}.
First we need to check a technical condition about the
$\aut(\p^N)$-action on $J_{g,Q}(X,d,t)$, namely that it is 
 {\it proper}. 
Properness of the group action is 
equivalent to the separatedness  of the
family considered in \cite[4.2]{FuPa}. 
We need the following:

\begin{prop}\label{separation.prop}
 In the setting of Theorem~\ref{modexists.general}
assume that $S$ is the spectrum of a complete DVR with quotient field
$K$ and residue field $k$. Let
$(C_K,f_K,\sigma_K,\rho_K)$ be a stable curve over $K$.
\begin{enumerate}
\item There is at most one way to extend it to a
stable curve $(C_S,f_S,\sigma_S,\rho_S)$ over $S$.
\item If $P\to S$ and $Q\to S$ are smooth, then
there is a scheme $S'$, which is also  the spectrum of a DVR,  and a local map 
$S'\to S$ such that the pull back of 
$(C_K,f_K,\sigma_K,\rho_K)$ can be extended to 
$(C_{S'},f_{S'},\sigma_{S'},\rho_{S'})$.
\end{enumerate}
\end{prop}

\noindent {\it Proof.}  For the first part, we plan to follow the proof of the
uniqueness of canonical models (see for instance
\cite[3.52]{ko-mo}).

If we have two extensions, they are dominated by a common
resolution. So it is enough to show how to recover
$C_S$ from a resolution of it. The notion of resolution here 
is slightly different
from the usual one since the general fiber may also be singular.
Thus we allow 3 types of operations:
\begin{enumerate}
\item blowing up a smooth point of the special fiber,
introducing a $(-1)$-curve,
\item blowing up an isolated singular  point of the special fiber,
introducing one or two  $(-2)$-curves, and
\item blowing up a nonisolated singular  point of the special fiber,
introducing  two  $(-1)$-curves in the normalization of $C_S$.
\end{enumerate}

We would like to consider the pair $(C_S, \sigma_S(P_S)+\rho_S(Q_S))$,
but it is not canonical. Instead we do the following.
Let $D_i$ be an irreducible component of $\sigma_S(P_S)+\rho_S(Q_S)$
with central fiber $(D_i)_k=\spec k'[t]/(t^{d_i})$
for some $k'\supset k$ and $d_i\geq 1$. (Here we use that  $S$ is complete.)
 Set $D=\sum \tfrac1{d_i}D_i$ as a $\q$-divisor
on $C_S$. The reason for introducing $D$ is  to assure
that every irreducible component of $D$ has local intersection number 1 with
$C_k$ at every geometric intersection point. 
Our stability assumption \ref{stable.gen.defn}.5
was made to assure that $\omega_{C_S/S}(D)$ is relatively ample.

We also see that under a blow up $C'_S\to C_S$, the exceptional
curves appear in $\omega_{C'_S/S}(D)$ with coefficient
0, unless we blow up a smooth point not on $D$, when 
the coefficient is 1. Thus if $g:C^m_S\to C_S$
is any resolution, the global sections of
 $(\omega_{C^m_S/S}(D))^{\otimes n}$ all pull back from
$(\omega_{C_S/S}(D))^{\otimes n}$ for $n\geq 1$. Thus
$C_S$ is the image of $C^m_S$ under the morphism given by the
 global sections of
 $(\omega_{C^m_S/S}(D))^{\otimes n}$ for some $n\gg 1$.

Assume next that $P\to S$ and $Q\to S$ are smooth.
In this case the existence proof of
\cite[4.2]{FuPa} works, except that one should use
the  semi stable reduction theorem in arbitrary characteristic
as proved in \cite{AW}.
\qed
\medskip

Thus we conclude  that ${\overline M}_{g,Q}(X,d,t)$
exists as an algebraic space that is
separated and of finite type  over $S$.
If $P\to S$ and $Q\to S$ are smooth,
the properness  of ${\overline M}_{g,Q}(X,d,t)$ over $S$ 
follows from Proposition~\ref{separation.prop}.  
The projectivity  over $S$ can be checked
as in \cite[3--4]{FuPa}, but it is actually not important
for any of what follows.

\end{say}

\begin{say}[Local structure of  ${\overline M}_{g,Q}(X,d,t)$]{\ }

We have constructed ${\overline M}_{g,Q}(X,d,t)$ as a 
quotient of $J_{g,Q}(X,d,t)$,
and so we expect to read off the local properties
 of ${\overline M}_{g,Q}(X,d,t)$ from the local properties of $J_{g,Q}(X,d,t)$.
The latter is actually much nicer, and  parts 1,2,3 and 5
of Theorem 
 \ref{modexists.prop.general}
follow from the corresponding properties of $J_{g,Q}(X,d,t)$.

\begin{prop}\label{J.local.prop}
 Let $S$ be a normal  scheme and
 $X$  a smooth projective scheme over $S$.
Let $P,Q$ be  finite  and flat  $S$-schemes 
of relative embedding dimension at most $1$ and let
$t:P\to X$ be a morphism. Let $k$ be a field,
$\spec k\to s\in S$ a point in $S$ and 
$$
(C,f:C\to X_k, \sigma:P_k\to C,\rho:Q_k\to C)
$$  a   stable curve 
of genus $g$ and degree $d$ with base point $t_k(P_k)$
and marking $Q_k$
over $k$. 
Let $L$ be the line bundle on $C$ defined as in 
Paragraph \ref{moduli.reduce.gen.rems}. Let 
$\delta:C\into \p^N_k$ be a nondegenerate embedding
by all sections of $L$
and  let $[(C,f,\sigma,\delta,\rho)]$ be the corresponding point in 
$J_{g,Q}(X,d,t)$.
 Assume  $H^1(C, L(-\sigma(P_k))=0$ and 
$H^1(C,f^*T_{X_k}(-\sigma(P_k))=0$.
Then
\begin{enumerate}
\item $ J_{g,Q}(X,d,t)\to \Hom^1(P,\p^N)$ is smooth at
 $[(C,f,\sigma,\delta,\rho)]$.
\item $ J_{g,Q}(X,d,t)$ has a unique irreducible component
$\comp(C,f,\sigma,\delta,\rho)$ 
passing through $[(C,f,\sigma,\delta,\rho)]$.
\item The dimension of $[(C,f,\sigma,\delta,\rho)]$ is
$$
-(f_*[C]\cdot K_X) +(\dim X-3)(1-g)+|Q| -|P|(\dim X-1)+\dim S+g+
(N+1)^2-1.
$$
\item A general point of $ \comp(C,f,\sigma,\delta,\rho)$
corresponds to a smooth stable curve.
\end{enumerate}
\end{prop}

\noindent {\it Proof.}
Recall from Claim~\ref{Q.gives.smooth.over} that
$ J_{g,Q}(X,d,t)\to W_{g}(X,d,t)$ is smooth of relative dimension
$|Q|$, thus the main point is to prove the analogous
statements for $ W_{g}(X,d,t)\to \Hom^1(P,\p^N)$.

The local structure of the Hilbert scheme is governed by the
normal bundle (see Paragraph   \ref{hilbert_schemes}). Thus the key 
to prove Proposition \ref{J.local.prop}
is to understand the normal bundle
of the image of 
$$
C\stackrel{(f,\delta)}{\longrightarrow}
 X\times_k \p^N\map B_{t'}(X\times_k \p^N).
$$
Going from the embedding into $ X\times_k \p^N$
to the blow up is straightforward, and the
important part is to describe the 
normal bundle
of the image of  $C$ in $ X\times_k \p^N$.

The argument  has three main ingredients.
\begin{enumerate}\setcounter{enumi}{3}
\item General remarks about the ``second exact sequence''.
\item Computing the normal bundle of a curve in $\p^N$.
\item Relating the normal bundle of $C\subset X\times \p^N$
to $f^*T_X$.
\end{enumerate}

\begin{say}[The second exact sequence]{\ }

Let $Z$ be a closed subscheme of $Y$ with ideal sheaf $\cI_Z$.
  The so-called ``second exact sequence''
(see \cite[II.8.4.A]{H77}) is
$$
 \cI_{Z}/\cI_{Z}^2\stackrel{\delta}{\to} \Omega^1_{Y}|_{Z}
\to \Omega^1_{Z}\to 0.
$$
In general it is not left exact, but it is easily seen to be left exact
at points where $Z$ and $Y$ are smooth.
If such points are dense in $Z$, $Z$ is reduced  and  $\cI_{Z}/\cI_{Z}^2$
is locally free, then $\delta$ can have  no kernel at all, so
 we also have left exactness. This holds if,
for instance, $Y$ is smooth and $Z$ is a reduced 
local complete intersection subscheme.

If $C\subset \p^N$ is a reduced 
local complete intersection curve, we can compute
the degree of its normal bundle from the second exact sequence
as
$$
\deg N_{C,\p^N}=(N+1)\deg C+2p_a(C)-2.
$$

\end{say}

\begin{say}[Normal bundles in $\p^N$]{\ }

Let $Z\subset \p^N$ be a reduced local complete intersection subscheme.
By the previous remarks, the natural map
$$
T_{\p^N}|_Z\to N_{Z,\p^N}
$$
is a surjection at smooth points of $Z$.
The tangent bundle of $\p^N$ is the quotient of 
a direct sum of $N+1$ copies of $\o_{\p^N}(1)$
(see \cite[II.8.13]{H77}), hence we get 
$$
\o_{\p^N}(1)^{N+1}|_Z\to N_{Z,\p^N}.
$$
If $Z=C$ is a curve, we get a surjection
of the first cohomologies
$$
H^1(C,\o_{\p^N}(1)|_C)^{N+1}\onto H^1(C,N_{C,\p^N}).
$$
Thus if $C$ is embedded by the sections of a line bundle $L$,
this allows us to get vanishing results about the normal bundle of $C$
in terms of  $H^1(C,L)$.
\end{say}

The normal bundle of a  diagonal embedding
is described by the following lemma.
We are interested in the case $Y=\p^N$.

\begin{lem}\label{diag.emb.seq}
 Let $X,Y$ be  smooth varieties over a field $k$.
Let $C$ be a proper, reduced, local complete intersection curve.
Let $f:C\to X$ be a morphism and $g:C\into Y$ an embedding.
Let $C'\subset X\times Y$ be the image of the diagonal embedding.
Then there is an exact sequence
$$
0\to f^*T_X\to N_{C',X\times Y}\to N_{C,Y}\to 0.
$$
\end{lem}

\noindent {\it Proof.} The middle term of the second exact sequence
$$
0\to \cI_{C'}/\cI_{C'}^2\to \Omega^1_{X\times Y}|_{C'}
\to \Omega^1_{C'}\to 0
$$
 can be identified as
$$
\Omega^1_{X\times Y}|_{C'}\cong f^*\Omega^1_{X}+ 
\Omega^1_{Y}|_{C}.
$$
Since $g$ is an embedding, the restriction 
$\Omega^1_{Y}|_{C}\to \Omega^1_{C}$ is surjective,
which implies that the composition
$$
\cI_{C'}/\cI_{C'}^2\to \Omega^1_{X\times Y}|_{C'}
\to f^*\Omega^1_{X}
$$
is surjective and its kernel is
the same as the kernel of $\Omega^1_{Y}|_{C}\to \Omega^1_{C}$,
which is just $\cI_C/\cI_C^2$ by the second exact sequence applied to
$C\subset Y$. Putting all these together we get the exact sequence
$$
0\to \cI_C/\cI_C^2 \to \cI_{C'}/\cI_{C'}^2 \to f^*\Omega^1_{X}\to 0,
$$
which is the dual of the claimed sequence.\qed

\end{say}

\begin{say}[The local structure of $J(X,d,t)$]{\ }

We can now put things together to
prove  parts 1,2 and 3 of Proposition~\ref{J.local.prop}.

Let $N_C$ be the normal bundle of
$C$ in $X_k\times\p^N$, embedded by $(f,\delta)$,
and let $N'_C$ be the   normal bundle of the birational transform of 
$C$ in $B_{t'}(X_k\times\p^N)$.
Then
 $N'_C\cong N_C(-\sigma_k(P_k))$.
Twisting the exact sequence of 
 Lemma~\ref{diag.emb.seq} we obtain  an exact sequence
$$
 0\to f^*T_X(-\sigma(P_k))\to N'_C\to N_{C,\p^N}(-\sigma(P_k))\to 0.
$$
We have proved that 
$N_{C,\p^N}(-\sigma(P_k))$ is generically a quotient of $N+1$ copies
of $L(-\sigma(P_k))$. This implies that
 $H^1(C,N'_C)=0$.
By Theorem~\ref{ksm}.2, this implies  
Proposition~\ref{J.local.prop}.1.

The scheme $\Hom^1(P,\p^N)$ is smooth over $S$. Therefore,
in a neighborhood of $\comp(C,f,\sigma,\delta,\rho)$, 
$W_{g}(X,d,t)\to \Hom^1(P,\p^N)$
is smooth over the normal scheme $S$,  hence itself normal.
A normal scheme is locally irreducible,
proving Proposition~\ref{J.local.prop}.2.

The dimension of the fiber of $W_{g}(X,d,t)\to \Hom^1(P,\p^N)$
at $[(C,f,\sigma,\delta,\rho)]$ is $H^0(C,N'_C)$. This is computed as
$$
\begin{array}{ll}
&\!\!\!\!\!\!\!\!\!  h^0(C,N'_C)= \chi(C,N'_C)\\
&= \deg N'_C+(\dim X+N-1)(1-g)\\
&= \deg f^*T_X(-\sigma(P_k))+\deg N_{C,\p^N}(-\sigma(P_k))+(\dim X+N-1)(1-g)\\
&= -(f_*[C]\cdot K_X)-|P|(\dim X+N-1)+
  (N+1)\deg L +2g-2+\\
&\hfill +(\dim X+N-1)(1-g)\\
&= -(f_*[C]\cdot K_X) +(\dim X-3)(1-g) -|P|(\dim X-1)+\\
& \hfill +(N+1)\deg L+N(1-g)-|P|\cdot N. 
\end{array}
$$

In order to get the dimension of $W_{g}(X,d,t)$,
we have to add the dimension of $\Hom^1(P,\p^N)$,
which is  $|P|\cdot N+\dim S$.
Finally using that $N=\deg L-g$ we get 
Proposition~\ref{J.local.prop}.3.
\end{say}

Proposition~\ref{J.local.prop}.4 is a special case of the next Lemma. \qed

\begin{lem} Let $Y$ be a smooth variety over $k$,
$C\subset Y$ a reduced curve with at worst nodes.
Assume that $N_{C,Y}$ is generated by global sections
and $H^1(C,N_{C,Y})=0$.

Then $\hilb(Y)$ has a unique irreducible component
containing $[C]$. This component is smooth at $[C]$ and
a nonempty open subset of it parametrizes
irreducible and smooth curves in $Y$.
\end{lem}

\noindent {\it Proof.} 
It is enough to prove that a general 1-parameter deformation
has no nodes near a given node of $C$.

This is a local question, so 
choose local (analytic or formal) coordinates so that $C$ is given by
$$
y_1y_2=y_3=\cdots=y_n=0.
$$
The coordinates $y_3,\dots,y_n$ generate a well-defined
subspace of $N_{C,Y}|_P$. Assume that our 1-parameter deformation
corresponds to a section of $N_{C,Y}$ that is not in that
subspace at $P$. Then we have equations
$$
y_1y_2+tg(t,{\mathbf y})=y_3+tg_3(t,{\mathbf y})=\cdots=
y_n+tg_n(t,{\mathbf y})=0,
$$
where $g(0,0)\neq 0$ by our assumption. 
We can choose new coordinates $y_i+tg_i(t,{\mathbf y})$
for $i=3,\dots,n$, and call them again $y_i$,
 to get  a simpler system
$$
y_1y_2+t(a+G(t,{\mathbf y}))=y_3=\cdots=
y_n=0,
$$
where $a\neq 0$ and $G(0,0)=0$.
The singular points are given by equations
$y_3=\cdots=y_n=0$ and
$$
y_2+t\dfrac{\partial G}{\partial y_1}=
y_1+t\dfrac{\partial G}{\partial y_2}=0.
$$
We can view these as equations
that can be substituted back into $y_1y_2+t(a+G(t,{\mathbf y}))=0$
to  get a new equation for
the supposed  singular point:
$$
ta=-tG(t,{\mathbf y})-
t^2\dfrac{\partial G}{\partial y_1}\dfrac{\partial G}{\partial y_2}.
$$
 The latter  has no solutions
for $t\neq 0$ and $t,y_1,\dots,y_n$ small since $a\neq 0$
and $G(0,0)=0$.\qed

\begin{say}[The case $g\geq 1$]

The only remaining problem is that we have
worked with $J_{g,Q}(X,d,t)$ but, as noted in \ref{J.const.say}.7,
 we would like
to conclude similar results about 
$J^*_{g,Q}(X,d,t)$.
The latter is a subscheme of $J_{g,Q}(X,d,t)$,
but it is better to view $J_{g,Q}(X,d,t)$
as $\pic^0$ of the universal curve over $J^*_{g,Q}(X,d,t)$.
By \cite{grot-pic},
$\pic^0$ of a flat family of curves is smooth,
hence we conclude that $J_{g,Q}(X,d,t)\to J^*_{g,Q}(X,d,t)$
is smooth of relative dimension $g$. Thus the smoothness results about
$J$ translate into smoothness results about $J^*$.

\end{say}

\begin{ack}   We   thank J.-L.\ Colliot-Th\'el\`ene and the referee
for a long list of corrections and suggestions.
Partial financial support was provided by  the NSF under grant number 
DMS-9970855 and DMS02-00883. The first author was partially supported by
CNPq.
\end{ack}


\bigskip

\noindent Mathematics Department,
Princeton University,
Princeton, NJ 08540.

\medskip 

\noindent caraujo@math.princeton.edu

\noindent  kollar@math.princeton.edu

\end{document}